\documentclass[a4paper,12pt,pdftex]{article}

\usepackage[applemac]{inputenc} 
\usepackage[T1]{fontenc} 
\usepackage{lmodern}
\usepackage[english]{babel}
\usepackage{xspace}
\usepackage{amsmath}
\usepackage{amssymb}
\usepackage[all,2cell,pdf]{xy} 
\UseAllTwocells
\SelectTips{cm}{11}
\CompileMatrices
\usepackage{xcolor}
\usepackage[amsmath,thmmarks,hyperref]{ntheorem}
\usepackage{calrsfs}
\usepackage{dsfont} 
\usepackage{mbboard}
\usepackage{mathtools}
\usepackage{extarrows}
\usepackage{chemarrow}
\usepackage{enumerate}
\usepackage{leftidx}
\usepackage{cmll}
\usepackage{verbatim}

\usepackage[pdftex]{hyperref}

\usepackage{cite}

\usepackage[textheight=26.7cm,hmargin=2.5cm]{geometry} 
\usepackage{cleveref}

\title{Parametrized $K$-Theory}
\author{Nicolas MICHEL\thanks{This research was started during my PhD studies at EPFL, Lausanne, Switzerland, and continued during my postdoc and lecturer positions at the same university. I was supported by the FNS grant 200020-121864.}}
\date{\small Section of mathematics\Par
EPFL\Par
Lausanne, Switzerland\\\bigskip
 \href{mailto:nicolas.michel@gmail.com}{nicolas.michel@gmail.com}
}

\theoremheaderfont{\medskip\normalfont\bfseries}
\theoremseparator{}
\theoremstyle{break}
\newtheorem{Prop}{Proposition}[section]
\newtheorem{Lem}[Prop]{Lemma}
\newtheorem{Thm}[Prop]{Theorem}
\newtheorem{Cor}[Prop]{Corollary}

\newtheorem{Axns}[Prop]{Axioms}

\theoremstyle{nonumberbreak}

\theoremstyle{plain}
\theoremseparator{ :}
\newtheorem{Def}[Prop]{Definition}
\newtheorem{Defp}[Prop]{Definitions}

\newtheorem{Not}[Prop]{Notation}

\theorembodyfont{\upshape}
\newtheorem{Rem}[Prop]{Remark}
\newtheorem{Remp}[Prop]{Remarks}

\newtheorem{Exp}[Prop]{Examples}

\theoremstyle{nonumberplain}
\theoremseparator{ :}

\theoremheaderfont{\slshape}
\theoremsymbol{\fbox{}}
\newtheorem{Pf}{Proof}
\qedsymbol{\PfSymbol}

\newenvironment{Exs}[1][]{%
\begin{Exp}#1\begin{enumerate}}%
{\end{enumerate}\end{Exp}}

\newenvironment{Defs}[1][]{%
\begin{Defp}#1
\begin{enumerate}}%
{\end{enumerate}\end{Defp}}

\newenvironment{Rems}[1][]{%
\begin{Remp}#1\begin{enumerate}}%
{\end{enumerate}\end{Remp}}

\newenvironment{Axioms}[1][]{%
\begin{Axns}#1\begin{enumerate}}%
{\end{enumerate}\end{Axns}}

\newcommand\N{\ensuremath{\mathds{N}}\xspace}
\newcommand\Z{\ensuremath{\mathds{Z}}\xspace}

\newcommand\R{\ensuremath{\mathds{R}}\xspace}

\newcommand{\Par}{\medskip\\}
\newcommand{\resp}{resp.\ }

\newcommand{\ra}{\xrightarrow}

\renewcommand{\epsilon}{\varepsilon}
\newcommand{\am}{\textnormal{(M)}\xspace}
\newcommand{\ac}{\textnormal{(C)}\xspace}
\newcommand{\al}{\textnormal{(L)}\xspace}
\newcommand{\aco}{\ensuremath{\mathrm{(C^{op})}}\xspace}
\newcommand{\alo}{\ensuremath{\mathrm{(L^{op})}}\xspace}
\newcommand{\tm}{\textnormal{($\widetilde{\mathrm{M}}$)}\xspace}
\newcommand{\tl}{\textnormal{($\widetilde{\mathrm{L}}$)}\xspace}
\newcommand{\tc}{\textnormal{($\widetilde{\mathrm{C}}$)}\xspace}
\newcommand{\cqfd}{\hfill\PfSymbol}

\newcommand{\C}{\ensuremath{\mathcal C}\xspace}
\newcommand{\D}{\ensuremath{\mathcal D}\xspace}
\newcommand{\E}{\ensuremath{\mathcal E}\xspace}
\newcommand{\B}{\ensuremath{\mathcal B}\xspace}
\newcommand{\T}{\ensuremath{\mathcal T}\xspace}
\newcommand{\V}{\ensuremath{\mathcal V}\xspace}
\newcommand{\A}{\ensuremath{\mathcal A}\xspace}
\newcommand{\one}{\ensuremath{\mathds 1}\xspace}
\newcommand{\two}{\ensuremath{\mathbb 2}\xspace}

\newcommand{\SYMMON}{\ensuremath{\mathit{SYMMON}}\xspace}

\newcommand{\Mod}{\ensuremath{\mathit{Mod}}\xspace}
\newcommand{\Mon}{\ensuremath{\mathit{Mon}}\xspace}
\newcommand{\Comm}{\ensuremath{\mathit{Comm}}\xspace}

\newcommand{\Free}{\ensuremath{\mathit{Free}}^{\mathit{fg}}\xspace}
\newcommand{\GBun}{\ensuremath{G\text{-}\mathit{Bun}}\xspace}
\newcommand{\VBun}{\ensuremath{\mathit{VBun}}\xspace}

\newcommand{\CAT}{\ensuremath{\mathit{CAT}}\xspace}

\newcommand{\F}{\ensuremath{\mathcal F}\xspace}

\newcommand{\Top}{\ensuremath{\mathit{Top}}\xspace}

\newcommand{\Triv}{\ensuremath{\mathit{Triv}}\xspace}

\newcommand{\Loc}{\ensuremath{\mathit{Loc}}\xspace}

\newcommand{\Sch}{\ensuremath{\mathit{Sch}}\xspace}
\newcommand{\Oc}[1]{\ensuremath{\mathcal{O}_{#1}}\xspace}

\newcommand{\Ab}{\ensuremath{\mathit{Ab}}\xspace}
\newcommand{\Ring}{\ensuremath{\mathit{Ring}}\xspace}
\newcommand{\Ringed}{\ensuremath{\mathit{Ringed}}\xspace}
\newcommand{\LRinged}{\ensuremath{\mathit{LRinged}}\xspace}
\newcommand{\Zar}{\ensuremath{\mathit{Zar}}\xspace}
\newcommand{\se}{\ensuremath{\mathrm{Ex}}\xspace}
\newcommand{\Aff}{\ensuremath{\mathit{Aff}}\xspace}

\newcommand{\G}{\ensuremath{\mathcal G}\xspace}
\newcommand\Ob{\ensuremath{\mathrm{Ob}\,}\xspace}
\newcommand\Mor{\ensuremath{\mathrm{Mor}\,}\xspace}
\newcommand\cod{\ensuremath{\mathrm{cod}}\xspace}
\newcommand\dom{\ensuremath{\mathrm{dom}}\xspace}
\newcommand{\Sh}{\ensuremath{\mathit{Sh}}\xspace}

\newcommand{\Spec}{\ensuremath{\mathrm{Spec}}\xspace}

\newcommand{\adjoint}[4]{\xymatrix{
	#3:#1\ar@<.9ex>[r]\ar@{}[r]|{\perp} & #2:#4 \ar@<.9ex>[l]
	}}

\begin{document}

\maketitle

\begin{abstract}
 In nature, one observes that a $K$-theory of an object is defined in two steps. First a ``structured'' category is associated to the object. Second, a $K$-theory machine is applied to the latter category to produce an infinite loop space. We develop a general framework that deals with the first step of this process. The $K$-theory of an object is defined via a category of ``locally trivial'' objects with respect to a pretopology. We study conditions ensuring an exact structure on such categories. We also consider morphisms in $K$-theory that such contexts naturally provide. We end by defining various $K$-theories of schemes and morphisms between them.
\bigskip
\\
\emph{Keywords}: $K$-theory -- Local triviality -- Exact categories -- Monoidal fibred categories -- Fibred Grothendieck sites -- Modules -- Sheaves of modules.\Par
2000 \emph{Mathematics Subject Classification}. Primary 18F25, 19D99; Secondary 13D15, 14F05, 14F20, 18D10, 18D30, 18D99, 18E10, 18F10, 19E08.

\end{abstract}

\tableofcontents

\section{Introduction}
\label{sec:introduction}

The definition of a $K$-theory of a mathematical object $X$ involves two steps. One first associates to $X$ a suitably structured category $\C_X$. One then applies to this structured category a $K$-theory machine that provides an infinite loop space $K(X)=K(\C_X)$.\footnote{Of course, there are also equivalent ad hoc constructions in some particular cases that do not follow this pattern. For instance, the $K$-theory of a ring can be obtained directly by the +-construction of Quillen and the $K$-theory of a space may be defined in a representable way \cite{WeiKbook}.} By ``structured category'', we mean a category with either a Quillen-exact (or more generally Waldhausen) structure, or a symmetric monoidal structure, which may be topologically enriched as well. For instance, the usual $K$-theories of a ring $R$ or a scheme $X$ are defined respectively via their exact categories of finitely generated projective $R$-modules and of locally free sheaves of $\Oc X$-modules of finite rank \cite{WeiKbook}. The $K$-theory of a space $X$ is obtained via its topologically enriched exact category of finitely generated projective $C(X)$-modules \cite{Pal96,Gra94,Gra05,Gil92,Mit01}. The usual $K$-theory of a ringed spectrum $R$ is defined via its Waldhausen category of finite cell $R$-modules \cite{EKMM97}.

In this article we take for granted the second step in this process, and concentrate on the first part, that is, the obtaining of a structured category from an object. We therefore naturally run into the following questions.
\begin{enumerate}
\item What kind of objects is $K$-theory designed to apply to?
\item Given such an object, what structured category should one associate to it in order to obtain a meaningful $K$-theory?
\item Finally, how does such a correspondence take the morphisms between these objects into account?
\end{enumerate}
After examining the available examples of $K$-theories in nature, we arrived at the following conclusion concerning the first question: the objects that typically admit a $K$-theory are commutative monoids, not only in a symmetric monoidal category, but in a symmetric monoidal \emph{opfibred} category $P$. More generally, they are objects of a category \B equipped with a functor
\begin{equation}
\label{eq:10}
F\colon\B\to\Comm(P)
\end{equation}
into a category of such monoids.

Let us turn to the second question. As we explain in \hyperref[sec:fibred-algebra]{section 2}, a symmetric monoidal opfibred category $P$ induces an opfibred category of modules over commutative monoids $\Mod_{c}(P)\to\Comm(P)$.\footnote{The notion of \emph{monoidal fibred category} is treated in \cite{Mal95,CA93,HM06,PS09,MyThesis}.} To an object $B$ of the category \B, we associate the category of modules over $F(B)$ (obtained by pulling back the latter opfibred category over the functor $F$ in \cref{eq:10}). The category of modules is not the category to which one wants to apply a $K$-theory functor, though, because it is unmanageable: it will in general not even be skeletally small. We obtain the desired categories of modules by restricting to ``locally trivial modules'', an idea inspired by the foundational article of Street \cite{Str04}.

Suppose the category \B is equipped with a pre-cotopology $J$, that is, basically, that each object $B$ of \B comes with a determined class of co-coverings (coverings on $\B^{op}$, see \cref{ssec:fibred-sites}). Suppose moreover that there are certain modules over certain monoids that are considered as trivial. An object $B\in\B$ is ``locally trivial'' if it can be covered in $J$ by trivial objects of \B. A module is ``locally trivial'' if its direct images over some $J$-covering are trivial modules. Given a locally trivial object $B\in\B$, we associate to it its full subcategory $\Loc_{B}$ of locally trivial modules. This is the right candidate for the $K$-theory machine, if it is skeletally small and if it has a structure of some kind that is accepted by a $K$-theory machine. In this article, we provide sufficient conditions to ensure the existence of an exact structure.

What about the third question? Under some hypotheses, there are not only categories of locally trivial modules over each locally trivial object $B\in\B$, but a \emph{sub-opfibration} of them. Therefore, a morphism of locally trivial objects in \B determines a functor between their categories of locally trivial modules. We study conditions ensuring that these functors are exact and thus induce morphisms in $K$-theory.

The theory we have developed to answer these questions happens to go beyond them. One can now modify existing $K$-theories by changing the diverse parameters in play as well as define $K$-theories in new areas. Moreover, under some natural hypotheses, the framework provides morphisms that allow us to compare the diverse $K$-theories obtained. It should also provide tools for studying descent questions.

This article is a summary and a continuation of the author's PhD thesis \cite{MyThesis}. For the sake of brevity, we do not mention set theoretical subtleties here, but an interested reader can have a look at \cite{MyThesis}. Moreover, this article assumes basic knowledge of the theory of (op-)fibred categories. We recommend the following references on this topic, if need be: \cite{BorII94,Shu08a,Gra69,Joh202,Str05,SGA1,JT94,Web07, MyThesis}. Note that we use interchangeably the synonyms ``(Grothendieck) fibration'' and ``fibred category''.

\begin{Not}\label{not:intro}
  In the same way that one sometimes considers a continuous map \linebreak $p\colon E\to B$ as a bundle with total space $E$ and base space $B$, one sometimes considers a functor ${P\colon\E\to\B}$ as a bundle of categories with \emph{total category} \E, \emph{base category} \B and \emph{fibre categories} the preimage categories $P^{-1}(B)$ for each $B\in B$. When we do so, we denote $P_{t}$, $P_{b}$ and $P_{B}$ the total, base and fibre categories of the functor $P$. The fibre of $P$ at $B\in\B$ is also denoted $\E_{B}$.
\end{Not}
\textbf{Acknowledgements} This article owes its existence to Prof.\ Kathryn Hess Bellwald, who, during my PhD studies at EPFL, has patiently guided me and encouraged me. I also thank her for reading several times this article and for her precious suggestions.

\section{Fibred algebra}
\label{sec:fibred-algebra}

In this section, we construct an opfibration of modules over monoids from a monoidal opfibration. One may skip this part on a first reading, and come back to it when needed.

We will see in the next section that many examples of opfibred sites with trivial objects arise as particular cases of this construction. We state the theory in the case of opfibrations rather than fibrations because the main examples require the opfibred setting. Moreover, the opfibred case is trickier, since extension of scalars is more complicated than restriction of scalars. The fibred case, which requires weaker hypotheses, is treated in \cite{MyThesis}.

\subsection{Monoidal opfibrations}
\label{ssec:mono-opfibr}

The notion of (strict) \emph{monoidal fibred category} is treated in \cite{Mal95}, whereas \emph{monoidal indexed categories} appear for instance in \cite{CA93,HM06,PS09}\footnote{Note that the \emph{monoidal fibrations} of \cite{Shu08a} are different notions, even though they have the same name.}. See \cite{MyThesis} for a systematic treatment of both the fibred and indexed frameworks and the description of a 2-equivalence between them. 
\begin{Def}
  A \emph{monoidal opfibration} over a category \B is a monoidal object in the 2-Cartesian 2-category of opfibrations over \B, opcartesian functors over \B and natural transformations over \B. It therefore consists in a sextuple $(\E\ra P\B,\otimes,I,\alpha,\lambda,\rho)$ where
  \begin{enumerate}
  \item $P\colon\E\to\B$ is an opfibration,
  \item
    $\xymatrix@!0@R=1.5cm@C=1cm{
      **[l]{\E \times_{\B}\E}\ar[dr] \ar[rr]^{\otimes} &&\E\ar[dl]^{P}\\
      &\B&
    }$

    is an opcartesian functor over \B,
  \item
    $\xymatrix@!0@R=1.5cm@C=1cm{
      \B\ar[dr]_{Id_{\B}} \ar[rr]^{I} &&\E\ar[dl]^{P}\\
      &\B&
    }$ 

    is an opcartesian functor over \B,
\item
\(\alpha\) is a natural isomorphism over \B, called the \emph{associator}, filling the following diagram
$$\xymatrixnocompile@!0@C=4pc@R=6pc{
&(\E \times_{\B}\E)\times_{\B}\E\ar[rr]^{\cong}\ar[ld]_{\otimes\times_{\B}Id_{\E}} &&
\E \times_{\B}(\E\times_{\B}\E)\ar[rd]^{Id_{\E}\times_{\B}\otimes} & \\
	\E\times_{\B}\E\ar[rrd]_{\otimes} &\rrtwocell<\omit>{!\mbox{$\stackrel{\alpha}{\cong}$}}&&&	\E\times_{\B}\E\ar[lld]^{\otimes}\\
&&\E&&
	}$$
\item
\(\lambda\) and \(\rho\) are natural isomorphisms over \B, called respectively the \emph{left\textnormal{ and }right unitors}, filling the following diagram.
$$\xymatrix@=1.5cm{
\B \times_{\B}\E\ar[r]^{u\times_{\B}Id_{\E}}\ar[dr]_{\cong}\drtwocell<\omit>{<-3>\lambda}&\E\times_{\B}\E\ar[d]^{\otimes} & \E\times_{\B}\B\ar[l]_{Id_{\E}\times_{\B}u}\ar[dl]^{\cong}\\
&\E\urtwocell<\omit>{<-3>\rho}&
}$$
\end{enumerate}
These data are subject to the usual \emph{coherence axioms} in each fibre of $P$.\\
Similarly, there is a notion of \emph{symmetric monoidal opfibration}.
\end{Def}
Given a monoidal opfibration $P$, each of its fibres is a monoidal category. Moreover, given a cocleavage on $P$, each direct image functor $f_{*}$ is a strong monoidal functor.

\subsection{Opfibred algebra}
\label{ssec:opfibred-algebra}

We now would like to do algebra in this setting. Recall that given a monoidal category \V, there is a fibration $\Mod(\V)\to\Mon(\V)$ of modules in \V over monoids in \V, whose inverse image functors are called \emph{restriction of scalars}. This fibration is a bifibration if \V has reflexive coequalizers and if the right tensors $-\otimes A$ preserve them (we consider right modules) \cite{Val03,MyThesis}. Its direct image functors are called \emph{extension of scalars}. We define modules in a monoidal opfibration directly in this fashion, because we are interested in the end in the opfibration of modules over monoids.

\begin{Defs}
  [Let $P\colon\E\to\B$ be a monoidal (resp.\ symmetric monoidal) opfibration.]
\item The \emph{category of monoids (resp.\ commutative monoids)} in $P$, denoted $\Mon(P)$ (resp.\ $\Comm(P)$) is specified as follows:
  \begin{itemize}
  \item Objects: (commutative) monoids in the fibres of $P$.
  \item Morphisms: A morphism $(R,\mu,\eta) \ra{\phi} (S,\nu,\lambda)$
    is a morphism $\phi\colon R\to S$ in \E, such that the two following diagrams commute.
$$\xymatrix@=1.3cm{
  R \otimes R \ar[r]^{\phi\otimes\phi}\ar[d]_\mu & S\otimes S\ar[d]^{\nu}\\
  R\ar[r]_\phi & S }\quad\quad \xymatrix@=1.3cm{
  I_{P(R)}\ar[d]_\eta \ar[r]^{u(P(\phi))} & I_{P(S)}\ar[d]^{\lambda}\\
  R\ar[r]_\phi & S }$$
\item Composition: composition of \E.
\end{itemize}
\item The \emph{category of modules} over monoids (resp.\ commutative monoids), denoted
\[\Mod(P) \mathit{\ (resp.\ } Mod_{c}(P)),\]
is specified as follows:
  \begin{itemize}
  \item Objects: Pairs $(R,M)$ where $R$ is a monoid (resp.\ commutative monoid) in $\E$ and $M$ is an $R$-module in $\E_{P(R)}$.
  \item Morphisms: Pairs $(R,(M,\kappa)) \ra{(\phi,\alpha)} (S,(N,\sigma))$ where $\phi\colon R\to S$ is in $\Mon(P)$ and $\alpha\colon M\to N$ is a morphism in \E such that: 
\begin{enumerate}
\item $P(\phi)=P(\alpha)$,
\item $\xymatrix@=1.3cm{
M\otimes R\ar[d]_\kappa\ar[r]^{\alpha\otimes\phi} & N\otimes S\ar[d]^\sigma\\
M \ar[r]_\alpha & N
}$
\end{enumerate}
\item Composition: Composition of $\E\times_{\B}\E$.
  \end{itemize}
\end{Defs}
One has thus a 2-story object $$\xymatrix{\Mod(P)\ar[d]\\\Mon(P)\ar[d]\\\B,}$$ where the top functor is the projection and the bottom functor is induced by $P$. For an object $B\in\B$, the restriction $\Mod(P)_{B}\to\Mon(P)_{B}$ is just the functor of modules over monoids $\Mod(\E_{B})\to\Mon(\E_{B})$.

\begin{Prop}\label{prop:modovermon}
Let $(\E\ra P\B,\otimes)$ be a monoidal opfibration. Suppose each fibre of $P$ has reflexive coequalizers and that they are preserved by direct image functors and tensors $-\otimes E$, for all $E\in \E$. Then, there is a sequence of opfibrations:
$$\xymatrix{\Mod(P)\ar[d]\\\Mon(P)\ar[d]\\\B.}$$
The same is true in the commutative setting.
\end{Prop}
\begin{Pf}
The functor $P$ induces an opfibration $Mon(P)\to\B$. Given a monoid $R$ over $A\in\B$ and an arrow $f\colon A\to B$ in \B, an opcartesian morphism $\underline{f}_{R}\colon R\to f_{*}R$ in $P$ is an opcartesian morphism in $\Mon(P)\to\B$ when $f_{*}R$ is given the unique monoid structure that makes $\underline f_{R}$ a morphism of monoids in $P$.

We construct now an opcartesian lift in $\Mod(P)\to\Mon(P)$ of a morphism of monoids $\phi\colon R\to S$ at a module $(R,M)$. Let $P(\phi)=f\colon A\to B$. The opcartesian lift 
\[{\underline f_{M}\colon M\to f_{*}M}\] 
of $f$ at $M$ in $P$ provides an opcartesian lift \[(\underline f_{R},\underline f_{M})\colon (R,M)\to(f_{*}R,f_{*}M)\] of $f$ at $(R,M)$ in the composite functor $\Mod(P)\to\B$, when $(f_{*}R,f_{*}M)$ is given the unique module structure that makes the pair $(\underline f_{R},\underline f_{M})$ a morphism of modules in $P$. 

Let us look now at the following diagram.
\[\xymatrix@C=1.5cm{
(R,M)\ar[dr]_{(\underline f_{R},\underline f_{M})}&(S,f_{*}M\otimes_{f_{*}R}S)& \Mod(P)\ar[dd]\\
&(f_{*}R,f_{*}M)\ar[u]_{(\underline\phi,\underline\phi_{\sharp})}&\\
R\ar[r]^{\phi}\ar[rd]_{\underline f_{R}} & S&\Mon(P)\ar[dd]\\
& f_{*}R\ar@{-->}[u]_{\underline\phi}&\\
A\ar[r]^{f} & B & \B
}\]
The morphism $(\underline\phi,\underline\phi_{\sharp})$ in $\Mod(P)$ is given by extension of scalars, which exists thanks to the hypotheses. It is opcartesian in the fibre functor $\Mod(\E_{B})\to\Mon(\E_{B})$. One can then prove, either directly or by using the dual of \cite[Lemma 2.1.45]{MyThesis}, that the following composite is opcartesian in $\Mod(P)\to\Mon(P)$:
\[(R,M)\ra{(\underline{f}_{R},\underline{f}_{M})}(f_{*}R,f_{*}M)\ra{(\phi,\phi_{\sharp})}(S,f_{*}M\otimes_{f_{*}R}S).\]
Finally, for the commutative case, one just has to check that commutative monoids are closed under direct images in $\Mon(P)\to\B$. 
\end{Pf}
\begin{Rems}
\item In the commutative setting, the opfibration of modules over commutative monoids is denoted by $\Mod_{c}(P)\to \Comm(P)$.
\item  The dual situation is easier because restriction of scalars is always defined. One can prove, in the same fashion, that given a monoidal \emph{fibration} $P$, one obtains a composite of fibrations $\Mod(P)\to\Mon(P)\to\B$ (this is proven in detail in \cite{MyThesis}). In particular, if a monoidal bifibration $P$ satisfies the hypotheses of \cref{prop:modovermon}, then the latter composite is a composite of bifibrations.
\end{Rems}

\begin{Exs}
\item If a monoidal category \V is considered as a monoidal opfibred category $\V\to\one$, then one recovers the usual construction of $\Mod(\V)\to\Mon(\V)$.
\item Given a category \C with pullbacks, the \emph{codomain functor} $\cod_{\C}\colon\C^{\two}\to\C$ from the arrow category of \C to \C is a monoidal fibration. It is also an opfibration, but not monoidal unless \C is a groupoid. If \C has all finite limits, then, given an internal group $G$, $G$-bundles, i.e., arrows $E\to B$ in \C with a fibrewise $G$-action, form a subcategory of $\Mod(\cod_{\C})$. See \cite{MyThesis} for more details.
\item There is an opfibration of sheaves of abelian groups over spaces $\Sh\to\Top^{op}$ that is obtained this way. Consider the \emph{fibration} given by the Grothendieck construction of the following indexed category over $\Top^{op}$:
\begin{equation}\label{eq:1}\begin{array}{rcl}
\Sh\colon(\Top^{op})^{op}	&\longrightarrow &CAT\\
X	& \longmapsto&\Sh(X)\\
f^{op}\colon Y\to X & \longmapsto& f_{*}\colon\Sh(X)\to\Sh(Y).
\end{array}
\end{equation}
Here $\Sh(X)$ is the category of sheaves of abelian groups over $X$ and $f_{*}$ is the direct image sheaf functor (which is the inverse image functor of this indexed category).
This fibration $\Sh\to\Top^{op}$ has the following form:
  \begin{enumerate}
  \item $\Ob\Sh$: Pairs $(X,\F)$ where $X$ is a space and $\F$ is sheaf of abelian groups on $X$.
  \item $\Mor\Sh$: A morphism $(f^{op},f^{\sharp})\colon (X,\F)\to(Y,\G)$ consists in a continuous map $f\colon Y\to X$ and a morphism $f^{\sharp}\colon \F\to f_{*}\G$ of sheaves of abelian groups over $X$.
  \end{enumerate}
It is also an opfibration by \cite[Th.\ 2.1.37]{MyThesis}, since for each continuous map $f$, there is an adjunction
\begin{equation}\label{eq:2}
\adjoint{\Sh(Y)}{\Sh(X)}{f^{-1}}{f_{*}}.
\end{equation}
In fact, this opfibration is symmetric monoidal, as we shall briefly explain now. The indexed category \eqref{eq:1} take values in the 2-category \SYMMON of symmetric monoidal categories. 
Its Grothendieck construction thus admits a structure of a symmetric monoidal object in $\CAT/\Top^{op}$  \cite{MyThesis}. It remains to show that its tensor product and unit are opcartesian functors, which follows from the fact that the oplax monoidal structure on $f^{-1}$ induced by the adjunction \eqref{eq:2} is opstrong \cite{MyThesis,LM09}.

The symmetric monoidal opfibration $\Sh\to\Top^{op}$ satisfies the conditions of \cref{prop:modovermon}. Indeed, there is a monoidal isomorphism $\Sh_{X}\cong\Sh(X)$. Now, the category $\Sh(X)$ is a cocomplete and closed monoidal category \cite{KS06}.
Moreover, the direct image functors are left adjoints because this opfibration is a bifibration. According to \cref{prop:modovermon}, there is thus a sequence of opfibrations, defined by modules over commutative monoids in the symmetric monoidal opfibration, which we write this way\footnote{We have chosen another notation here than in \cite{MyThesis}: what is called $\Oc{}\text{-}\Mod$ here is its dual in \cite{MyThesis}.}:
\begin{equation}\label{eq:3}
\Oc{}\text{-}\Mod\to\Ringed^{op}\to\Top^{op}.
\end{equation}
The category \Ringed is precisely the category of \emph{ringed spaces} of the literature, e.g. \cite{Qin02,Har77}. For each ringed space $(X,\Oc X)$, the fibre $\Oc{}$-$\Mod_{(X,\Oc X)}$ is the category of \emph{sheaves of \Oc X-modules} of the literature. In the opfibration $\Oc{}\text{-}\Mod\to\Ringed^{op}$, the direct image of an object $(Y,\Oc Y,\G)$ over a morphism $(f,f^{\sharp})\colon(X,\Oc X)\to(Y,\Oc Y)$ in \Ringed is the object $(X,\Oc X,f^{-1}\G\otimes_{f^{-1}\Oc Y}\Oc X)$\footnote{It is called the \emph{inverse image} of $\G$ in the literature. Indeed, in the literature one considers an indexed category over \Ringed that does not come from the dual of this opfibration (the dual would have the opposite categories of modules). We think it is natural to consider the opfibration $\Oc{}\text{-}\Mod\to\Ringed^{op}$, rather than the usual fibration over \Ringed, even though one has been used to work with the category of ringed spaces (versus its dual). Indeed, the so-called inverse image along a morphism $(f,f^{\sharp})$ in \Ringed is really a direct image, since it is a left adjoint functor between the categories of modules. Note moreover that it involves an extension of scalars, not a restriction. In addition, the inverse image $f^{-1}\G$ of a sheaf $\G$ over $Y$ with respect to a continuous map $f\colon X\to Y$ is fundamentally the direct image with respect to the functor $f^{-1}\colon\T(Y)\to\T(X)$ between the topologies. It is indeed a left adjoint functor between the categories of sheaves.}. 

\end{Exs}

\section{Locally trivial objects}
\label{sec:locally-triv-objects}

In this section, we describe a general procedure for obtaining ``locally trivial objects'' in a ``fibred site with trivial objects''. The basic idea is the following. One has a site \B and a Grothendieck fibration $\E\to\B$ (or equivalently a pseudofunctor $\B^{op}\to\CAT$) over this site. For instance, this fibration can be the dual of the opfibration $\Mod_{c}(P)\to\Comm(P)$ defined in \cref{ssec:opfibred-algebra}. Now, some objects $B$ in the site \B and some objects in their fibre category $\E_{B}$ are called ``trivial''. In practice, these objects have specific properties that make them simpler to deal with. We want to study objects $B$ that are ``locally'' trivial in the site \B via their associated category $\Loc_{B}$ of ``locally trivial objects'' in $\E_{B}$. More precisely, we are looking for conditions insuring an exact structure on such categories in order to apply to them a $K$-theory machine. A locally trivial object $B$ in the site is thus studied via its $K$-theory defined by $K(B)=K(\Loc_{B})$.

We present the theory in the fibred context rather that the opfibred one. Indeed, the fibred theory is more intuitive because it is the context of the geometric examples such as vector bundles or principal bundles. Moreover, it allows us to work with sites, a much more common notion than cosites, and thus to avoid a proliferation of ``co'' and ``op'' that would obscure the text. Nevertheless, as we have already mentioned, important examples are expressed in the opfibred setting. When dealing with such situations, one has two choices. First, one can dualize the theorems, and the hints given throughout the text should facilitate this fairly straightforward process. Second, if one prefers to work with fibrations only, one can consider the dual fibration $P^{op}$ of the opfibration $P$ of study. In this case, one must be aware that one obtains results on the dual of the categories of locally trivial objects that the opfibration would have produced. When working in the realm of $K$-theory of exact categories, this does not really matter though, since the dual of an exact category \C is also exact and the $Q$-construction does not differentiate \C from its dual \cite[Ch. 4]{WeiKbook}.

Our notion of locally trivial object is inspired by that of Street in the foundational article \cite{Str04}. Our locally trivial objects are more general though, and used in a quite different perspective\footnote{See \cite{MyThesis} for a precise comparison with Street's notion of local triviality.}.

\subsection{Trivial objects}
\label{ssec:trivial-objects}
We start by defining what we mean by trivial objects in a fibration. Recall that a functor $F\colon\A\to\B$ is \emph{replete} if, given an isomorphism $g\colon F(A) \ra{\cong}B$ in \B, there exists an isomorphism $f\colon A\ra\cong A'$ in \A such that $g=F(f)$. A subcategory $\A\subset \B$ is \emph{replete} if its inclusion functor is so, that is, if given an object $A\in\A$ and an isomorphism $f\colon A\ra{\simeq}B$ in \B, both $B$ and $f$ belong to \A. Moreover, the \emph{replete full image} of a functor $F\colon\A\to\B$ is the full subcategory of \B of all objects isomorphic to some $F(A)$, $A\in\A$. 
\begin{Defs}
[\label{def:triv} Let $P\colon\E\to\B$ be a fibration.]
\item   A \emph{subfunctor of trivial objects} of $P$ is a subfunctor
$$\xymatrix@=3em{
\Triv_{t}\ar[d]_{\Triv}\ar@{^{(}->}[r] & \E\ar[d]^{P}\\
\Triv_{b}\ar@{^{(}->}[r] & \B
}$$ 
of $P$ that is:
  \begin{itemize}
  \item \emph{globally full}: $\Triv_{t}\subset \E$ and $\Triv_{b}\subset\B$ are full subcategories,
  \item \emph{replete}: Each fibre of \Triv is a replete subcategory of the corresponding fibre of $P$.
  \end{itemize}
Objects of both the total and base categories $\Triv_{t}$ and $\Triv_{b}$ are called \emph{trivial}.
\item\label{def:subfib} A \emph{subfibration} of $P$ is a subfunctor $Q\subset P$ that is a fibration and whose inclusion functor is Cartesian. It is a \emph{subfibration over \B} if $Q_{b}=\B$ (recall \cref{not:intro}).
\item \label{def:subfibtriv} A \emph{subfibration of trivial objects} is a subfunctor of trivial objects that is a subfibration. Equivalently, it is a globally full subfunctor $\Triv\subset P$ such that objects of $\Triv_{t}$ are closed under taking inverse images in $P$ over arrows in $\Triv_{b}$.
\end{Defs}

In important cases, the subfunctor of trivial objects arises as some sort of replete full image of a morphism of fibrations (morphisms of fibrations are not supposed to be Cartesian in this article):
\[\xymatrix{
\D\ar[r]^{F}\ar[d]_{Q} & \E\ar[d]^{P}\\
\A\ar[r]_{G} & \B.
}\]
There are two versions: the \emph{global replete full image} and the \emph{replete full image} of $(F,G)$. The former consists in the restriction of $P$ to the replete full images of both $F$ and $G$. The latter consists in the restriction of the domain of $P$ to the full subcategory of objects \textit{vertically} isomorphic to an $F(D)$ and in the restriction of its codomain to the full image of $G$. These two concepts coincide when $G$ is a replete functor. Note that these images of a morphism of fibrations are subfunctors but not necessarily subfibrations. They do produce a subfibration when $(F,G)$ is Cartesian and $G$ full \cite{MyThesis}. In a number of examples, the base functor $G$ is just the identity (and consequently, the two notions of image introduced coincide). 

There is a typical situation that automatically produces a subfibration of trivial objects over \B: when the fibration $P\colon\E\to\B$ has the property that its base category has a terminal object $*$. Indeed, given a class $\T$ of objects in the fibre $\E_{*}$ over the terminal object $*$, one defines trivial objects as follow. An object $E\in\E$ is trivial if there exists a Cartesian arrow $E\to T$ over $P(E)\to *$ with $T\in\T$.

We now give some examples of trivial objects in a fibration, most of which are determined by a subclass of the fibre over a terminal object.
\begin{Exs}\label{ex:triv}
\item (\textit{Trivial bundles}) Let \C be a category with finite limits. The codomain functor $\cod\colon\C^{\two}\to\C$ is a fibration. Given some subclass $\T$ of $\Ob\C$, trivial bundles are defined as above via the fibre over the terminal object (${\C/*\cong\C}$). These are precisely the product bundles (with fibre in \T). We next consider important algebraic variants.
\item\label{ex:triv:gbun} (\textit{Trivial $G$-bundles}) For any group object $G$ in a category \C with finite limits, there is a fibration $\GBun(\C)\to\C$ of $G$-bundles. Trivial objects (product $G$-bundles) are obtained by considering the subclass $\T=\{G\}$ of the $G$-objects in \C.

\item\label{ex:triv:vbun} (\textit{Trivial bundles of vector spaces}) There is a fibration $\VBun\to\Top$ of bundles of real vector spaces. Trivial bundles (product bundles) are obtained via the subclass $\{\R^{n}\mid n\in\N\}$ of topological vector spaces. The complex case is similar.

\item\label{ex:triv:csheaves} (\textit{Constant sheaves}) The fibre over a terminal object of the fibration of ringed spaces $\Ringed\to\Top$ defined in \cref{eq:3} is isomorphic to $\Comm^{op}$. The trivial objects determined by the whole fibre $\T=\Comm^{op}$ are \emph{constant sheaves}.

\item\label{ex:triv:affine} (\textit{Affine schemes}) Recall that a ringed space $(X,\Oc X)$ is \emph{locally ringed} if each stalk $\Oc {X,x}$, $x\in X$, is a local ring. Moreover, given two locally ringed spaces $(X,\Oc X)$ and $(Y,\Oc Y)$, a morphism of ringed spaces $(f,f^\sharp)$ between them is a \emph{morphism of locally ringed spaces} if the induced homomorphism of rings $$(f_\sharp)_x\colon\Oc{Y,f(x)}\to\Oc{X,x}$$ is local \cite{Har77,Qin02,MyThesis}. The fibration of ringed spaces $\Ringed\to\Top$ defined in \cref{eq:3} restrict to the fibration $\LRinged\to\Top$ where \LRinged is the category of locally ringed spaces. Now, there is a (non-Cartesian) morphism of fibrations into the fibration of locally ringed spaces given by:
\[\xymatrix@=1.3cm{
\Comm^{op} \ar[r]^{(\mathrm{Spec},\Oc{})} \ar[d]_{Id_{\Comm^{op}}}& \LRinged\ar[d]\\
\Comm^{op} \ar[r]_{\Spec} & \Top.
}\]
The global replete full image of this morphism determines the subfunctor of \emph{affine schemes}. We denote \Aff the category of affine schemes.

\item\label{ex:triv:freemod} (\textit{Free modules over rings}) Let us consider the opfibration of modules over rings $\Mod\to\Ring$ (we work in the dual framework of an opfibred category in this example). The sub-opfibration of trivial objects of interest here is that of finitely generated free modules. It is obtained via the subclass $\{\Z^{n}\mid n\in\N\}$ of $\Ob\Mod_{\Z}$, \Z being an initial object of \Ring. In the sequel, we will consider the restriction of this example to commutative rings.

\item\label{ex:triv:fp} (\textit{Finitely presented modules over rings}) Finitely presented modules form a sub-opfibration of trivial objects in the opfibration $\Mod\to\Ring$. Indeed, extension of scalars functors are right exact since they are left adjoints.

\item\label{ex:triv:freemodringed} (\textit{Free modules over ringed spaces}) We consider here the opfibration of sheaves of modules over ringed spaces $\Oc{}$-$Mod\to\Ringed^{op}$ as defined in \cref{eq:3}. A terminal object in \Ringed is given by $(*,\Z)$. When choosing $\T=\{(*,\Z,\Z^{n})\mid n\in\N\}$, the trivial objects are the \emph{free modules of finite rank}, that is, $\Oc X$-modules that are finite sums of $\Oc X$.
\item\label{ex:triv:affshemod} (\textit{Affine sheaves of modules over schemes}) The opfibration is the restriction 
\[\Oc{}\text{-}Mod_{l}\to\LRinged^{op}\] 
of the previous opfibration to the (non-full) subcategory $\LRinged\subset\Ringed$ of locally ringed spaces and their morphisms. Trivial objects are determined by the following morphism of opfibrations:
\begin{equation}\label{eq:4}\begin{aligned}\xymatrix@=1.3cm{
\Mod_{c} \ar[r]^-{(\mathrm{Spec},\Oc{},\tilde{\ })} \ar[d]& \Oc{}\text{-}\Mod_l\ar[d]\\
\Comm \ar[r]_-{(\mathrm{Spec},\Oc{})} & \LRinged^{op},
}\end{aligned}\end{equation}
where the category $\Mod_{c}$ is the restriction of $\Mod$ over \Comm and, for any $A$-module $M$, the sheaf $\tilde{M}$ is the associated sheaf of $\Spec A$-modules. The global replete full image of this morphism is a sub-opfibration of trivial objects because this morphism is opcartesian and has a full base functor \cite{MyThesis}. Trivial objects in the base are precisely affine schemes. We call the trivial modules \emph{affine modules}. In some applications, one restricts the domain opfibration in \cref{eq:4} to finitely presented modules or finitely generated free modules. In the first case, we call the trivial modules \emph{finitely presented affine modules}. In the second case, one gets free modules of finite rank as described in the previous example.

\end{Exs}

\subsection{Fibred sites}
\label{ssec:fibred-sites}

In order to introduce the notion of locally trivial objects, we must now define the notion of fibred site. Let us first state some basic definitions in order to fix the terminology and notation.

\begin{Defs}
\item A \emph{covering} of an object $C$ in a category \C is a set of arrows with codomain $C$.
\item A covering $S$ \emph{refines} a covering $R$ if each arrow of $S$ factors through some arrow of $R$.
\item A \emph{site} is a category \C (not necessarily skeletally small) together with a \emph{covering function}, i.e., a function $J$ that associates to each object $C$ of \C a class $J(C)$ of coverings of $C$. These coverings are called $J$-coverings.
\item A \emph{morphism of sites} $F\colon(\B,I)\to (\C,J)$ is a functor $F\colon \B\to\C$ such that for each $I$-covering $R$, the covering $F(R)$ is refined by a $J$-covering.
\item Let $J$ and $J'$ be covering functions on a category \C. $J$ is \emph{subordinated} to $J'$, and we write $J\preceq J'$, if the identity functor is a morphism of sites $Id\colon(\C,J)\to(\C,J')$.
\item Covering functions $J$ and $J'$ are \emph{equivalent} if both $J\preceq J'$ and $J'\preceq J$ hold.
\end{Defs}

We will require from place to place various axioms on covering functions. For each axiom, we give first a weak version that is good enough most of the time, and then a stronger version that is occasionally needed.

\begin{Axioms}
[Let $J$ be a covering function on \C.] 
\item[\normalfont\am]\emph{(Maximality)} For each $C\in \C$, $\{1_{C}\}\in J(C)$. 
\item[\tm]\emph{(Maximality, stronger version)} All isomorphisms belong to $J$.
\item[\normalfont\ac]\emph{(Coverage)} Given a $J$-covering $R$ of an object $C$ and an arrow $f\colon B\to C$, there exists a $J$-covering $S$ of $B$ such that the composite covering      \begin{equation}\label{eq:5}
    f\circ S\coloneqq\{f\circ g\mid g\in S\}
  \end{equation}
 refines $R$.
\item[\tc]\emph{ (Coverage, stronger version)} For every $J$-covering $R$, every $C_f\ra{f}C\in R$ and every arrow $g\colon D\to C$, there exists a pullback
\begin{equation}\label{eq:6}
\begin{aligned}
\xymatrix@!{
g^*(C_f) \ar[d]_{\bar f} \ar[r] & C_f \ar[d]^{f}\\
D \ar[r]_g & C
}\end{aligned}\end{equation}
such that the covering $\{g^*(C_f)\ra{\bar f} D\mid f\in R\}$ belongs to $J$.
\item[\normalfont\al] \emph{(Local character)} Given a $J$-covering $R$ and, for each $f\in R$, a $J$-covering $R_{f}$ of $\dom f$, the composite covering $\bigcup_{f\in R}f\circ R_{f}$ is refined by a $J$-covering.
\item[\tl] \emph{(Local character, stronger version)}
Given a $J$-covering $R$ and, for each $f\in R$, a $J$-covering $R_f$ of $\dom f$, the composite covering $\bigcup_{f\in R}f\circ R_{f}$ belongs to $J$.
\end{Axioms}

\begin{Defs}\label{def:pretop}
  \item A covering function satisfying axioms \ac is called a \emph{coverage}.
  \item We call \emph{pretopology} a covering function satisfying the three axioms \am, \ac and \al.
  \item It is a \emph{Grothendieck pretopology} if it satisfies their stronger versions\footnote{When the category \C has pullbacks, Grothendieck pretopologies are not really stronger than pretopologies. Indeed, the ``saturation'' of a pretopology $J$, which is equivalent to $J$, is a Grothendieck pretopology. See \cite{MyThesis} (beware of some changes in terminology though).}.
\end{Defs}

\begin{Rem}[Dualization]
 The opfibred setting requires the dual notions. All these definitions can be easily dualized, bringing in particular
  the notions of \emph{co-covering}, \emph{co-covering function},
  \emph{cosite} and \emph{pre-cotopology}; the ``co-notion'' on \C is
  to be defined so that it gives rise to the notion on $\C^{op}$. To each axiom of covering functions corresponds a dual version for co-covering functions. Except for \am which is self dual, we write \aco, \alo, …
\end{Rem}

We now describe the pretopologies we are going to use in the sequel.

\begin{Exs}
[\label{ex:pretop}Let \C be a category.]
\item \textit{The coarsest pretopology on \C:} This pretopology is denoted $\mathit{Coarsest}$ and has, for each object $C\in\C$, the identity-covering $\{1_{C}\}$.
\item \textit{The finest pretopology on \C}: What one usually calls the finest pretopology on \C is the covering function that has all possible coverings of objects of \C. It happens to be equivalent to a much simpler one, the covering function that has, for each object $C\in\C$, the identity covering and the empty covering $\emptyset$\footnote{The former is indeed the \emph{saturation} of the latter \cite{MyThesis,Vis08}.}. This is what we call the the finest pretopology and denote $\mathit{Finest}$.

Note that for every covering function $J$ satisfying axiom $\am$, one has the following relation
\[\mathit{Coarsest}\preceq J\preceq\mathit{Finest}.\]
\item \textit{The Zariski Grothendieck pre-cotopology on the category of commutative rings:} Its coverings are indexed sets $\{R\to R[a_{i}^{-1}]\}_{i\in I}$, where the set $\{a_{i}\}_{i\in I}$ generates $R$ as an ideal ($R=(a_{i})_{i\in I}$) and each $R\to R[a_{i}^{-1}]$ is a localization of $R$ at $a_{i}\in R$. It is denoted \Zar.
\item \textit{The open subset pretopology on the category of topological spaces:} It has, for each space $X$, all families $\{U_{i}\}_{i\in I}$ of open subsets of $X$ such that $X=\bigcup_{i\in I}U_{i}$.
\item\label{ex:pretop:sch} \textit{The Zariski pretopology on the category of ringed spaces:} It has, for each ringed space $(X,\Oc X)$, all families $\{(U_{i},\Oc X|_{U_{i}}\}_{i\in I}$ such that $\{U_{i}\}_{i\in I}$ is a covering of $X$ in the open subset pretopology. It is denoted \Zar. This pretopology is equivalent to the Grothendieck pretopology of families of jointly surjective open immersions. 
\item \textit{Pretopologies on the category of schemes:} First, the pretopology \Zar on \Ringed restricts to a pretopology on the category \Sch of schemes since an open subset of a scheme is a scheme \cite{Qin02}. In addition, we define the following Grothendieck pretopologies on \Sch \cite[Ch.\ 30]{StaProj}:
  \begin{itemize}
  \item \textit{étale}: An étale covering is a jointly surjective covering by étale morphisms.
  \item \textit{smooth}: A smooth covering is a jointly surjective covering by smooth morphisms.
  \item \textit{syntomic}: A syntomic covering is a jointly surjective covering by syntomic morphisms.
  \item \textit{fpqc}: An \textit{fpqc}-covering of a scheme $(X,\Oc X)$ is a covering \[\{(f_{i},f_{i}^{\sharp})\colon (X_{i},\Oc{X_{i}})\to(X,\Oc X)\}_{i\in I}\] by flat morphisms such that, for every open affine subscheme $(U,\Oc X|U)$, there are $i_{1},…,i_{n}\in I$ and open affine subschemes $(V_{j},\Oc{X_{i_{j}}}|V_{j})$,  $1≤j≤n$, with \[\bigcup_{j=1}^{n}f_{i_{j}}(V_{j})=U.\]
  \item \emph{fppf}: An \emph{fppf}-covering is a jointly surjective covering of flat morphisms locally of finite presentation.
  \end{itemize}
One has $\mathit{Zar}\subset\mathit{\acute etale}\subset\mathit{smooth}\subset\mathit{syntomic}\subset\mathit{fppf}\subset\mathit{fpqc}$ \cite{StaProj}.
\end{Exs}

\begin{Def}\label{def:fibsite}
A \emph{fibred site} $P\colon\E\to(\B,J)$ is a fibred category $P\colon\E\to\B$ together with a covering function $J$ on its base.
\end{Def}

\begin{Rem}\label{rem:fibsite}
 The notion of fibred site has already appeared in literature, for instance in \cite{Jar06,SGAIV2}. Our notion coincides with that of \cite{Jar06}, except that its author works in a ``sifted'' context. On the other side, the fibred sites of \cite{SGAIV2} are quite different mathematical structures. The latter paper deals with objects that could be described as ``categories fibred in sites'', whereas our notion is closely related to internal fibrations in some 2-category of sites. Indeed, the covering function on the base category \B of a fibred site $P\colon\E\to(\B,J)$ induces a covering function $J_{P}$ on \E by considering all Cartesian lifts of $J$-coverings. When \E is equipped with this covering function, the fibration $P$ becomes an internal fibration in some 2-category of sites. See \cite{MyThesis} for more details.
\end{Rem}

\subsection{Locally trivial objects}
\label{ssec:locally-triv-objects}

We are now ready to define locally trivial objects.

\begin{Def}
Let $\Triv\subset P\colon\E\to\B$ be a subfunctor of trivial objects.
\begin{itemize}
\item An object $B\in\B$ is \emph{locally trivial (in the base category) for a covering $R$} if $R$ is a covering of $B$ by trivial objects.
\item An object $E\in\E$ is\emph{ locally trivial (in the total category) for a covering $R$} if $R$ is a covering of $P(E)$ such that the inverse images of $E$ along the arrows of $R$ are trivial.
\end{itemize}
\end{Def}

The following elementary lemma will be of great use in the sequel.
\begin{Lem}\label{lem:refinetriv}
  Let $\Triv\subset P$ be a subfibration of trivial objects. If $E$ is locally trivial for a covering $R$, it is so for any of its refinements having trivial domains.\cqfd
\end{Lem}

\begin{Defs}
[Let $(P,J)$ be a fibred site such that $J$ satisfies axiom \am and \Triv a subfunctor of trivial objects.]
\item An object $B\in\B$, \resp $E\in\E$, is \emph{locally trivial (in $J$)} if it is locally trivial for some $J$-covering. 

The full subcategories of \B and \E consisting of these objects provide us with a subfunctor $\Loc\subset P$ of locally trivial objects:
$$\xymatrix@=3em{
\Triv_{t}\ar[d]_{\Triv}\ar@{^{(}->}[r] & \Loc_{t}\ar[d]_{\Loc}\ar@{^{(}->}[r] &\E\ar[d]^{P}\\
\Triv_{b}\ar@{^{(}->}[r] & \Loc_{b}\ar@{^{(}->}[r] &\B
}$$
\item The covering function $J$ on \B induces a covering function $J_{l}$ (``l'' for ``local'') on $\Loc_{b}$ whose coverings are $J$-coverings by trivial objects. It restricts to $\Triv_{b}$ and also extends to the whole category \B by giving no covering to objects that are not locally trivial. In order to avoid confusion, we always specify the site's category in these cases: $(\Triv_{b},J_{l})$ and $(\B,J_{l})$.
\end{Defs}

One has thus a sequence of subfunctors $\Triv\subset\Loc\subset P$. The subfunctor of locally trivial objects vary from \Triv to $P$ when $J$ varies from the coarsest to the finest pretopology  on \B. Note that locally trivial objects in the total category admit a characterization similar to those in the base category: they are precisely the objects that can be $J_{P}$-covered by trivial objects (c.f. \cref{rem:fibsite}).

As we shall see, most properties of locally trivial objects rely on axioms on the covering function $J_{l}$ rather than $J$. Yet, there are common situations where the properties \ac and \al pass  from $J$ to $J_{l}$. Of course, this is the case when $\Triv_{b}=\B$ (then $J=J_{l}$) and if one is interested only in these types of examples, which quite often arise, one can just forget about the index $l$ in the following results. Yet there are more subtle situations where $J≠J_{l}$, but where axioms \ac and \al pass  from $J$ to $J_{l}$. This is for instance the case for schemes. Here is a criterium for this to happen.
\begin{Lem}\label{lem:J_l}
  \begin{enumerate}[(i)]
  \item \label{lem:J_l:i} Suppose every $J$-covering of \emph{locally trivial} objects in the base admits a $J_{l}$-refinement. Then if $J$ satisfies axiom \ac (\resp \al), so does $J_{l}$.
  \item \label{lem:J_l:ii} Suppose that $J$ satisfies axioms \am, \ac and \tl and that every $J$-covering of \emph{trivial} objects in the base admits a $J_{l}$-refinement. Then $J_{l}$ satisfies axioms \ac and \al.
  \end{enumerate}
\end{Lem}

\begin{Pf}
  The first assertion is easy to verify. For the second statement, we first show that under these hypotheses, every $J$-covering of a \emph{locally} trivial object admits a $J_{l}$-refinement. We can then apply the first statement.\\
Let $B\in\Loc_b$ and $R$ a $J$-covering of $B$. Since $B$ is locally trivial, it admits a $J_{l}$-covering $S$. Let $g\colon \tilde B\to B$ be an arrow of this covering $S$. By axiom \ac, there exists a $J$-covering $S'_g$ of $\tilde B$ such that the composite covering $g\circ S'_g$ refines $R$. Moreover, by definition of a $J_{l}$-covering, ${\tilde B\in\Triv_b}$. Therefore, by hypothesis, there exists a $J_{l}$-covering $S_g$ of $\tilde B$ that refines $S'_g$. Thus, $\bigcup_{g\in S}(g\circ S_g)$ refines $\bigcup_{g\in S}(g\circ S'_g)$, which refines $R$. Consequently, $\bigcup_{g\in S}(g\circ S_g)$ refines $R$. Furthermore, since $J$ satisfies ($\tilde{\mathrm L}$), $\bigcup_{g\in S}(g\circ S_g)$ is a $J$-covering, and thus a $J_{l}$-covering.
\end{Pf}

\begin{Def}\label{def:fibsitetriv}
  A \emph{fibred site with trivial objects} $(P,J,\Triv)$ is a \hyperref[def:fibsite]{fibred site} $(P,J)$ equipped with a \hyperref[def:subfibtriv]{subfibration of trivial objects} \Triv, such that:
  \begin{itemize}
  \item $J$ satisfies \am,
  \item $J_{l}$ satisfies \ac.
  \end{itemize}
\end{Def}

The following lemma will be very useful later on.
\begin{Lem}\label{lem:comtriv}
  Let $(P,J,\Triv)$ be a fibred site with trivial objects such that $J_{l}$ satisfies \al. Let $E$ and $E'$ be locally trivial objects over the same object $B$.\Par
  Then, there exists a $J$-covering that trivializes both $E$ and $E'$.
\end{Lem}
\begin{Pf}
  Since $J_{l}$ is supposed to satisfy both axioms (C) and (L), this follows by lemma 2.2.8 in \cite{MyThesis} and \cref{lem:refinetriv}.
\end{Pf}

We come now to the main result of this section.

\begin{Prop}\label{prop:loc}
  Let $(P,J,\Triv)$ be a fibred site with trivial objects. Then $\Loc\subset P$ is a subfibration of $P$. Moreover, it is replete and globally full.
\end{Prop}

\begin{Pf}
Let $E$ be locally trivial for a covering $R$ of $B=P(E)$, and $g\colon B'\to B$ a morphism in $\Loc_b$. Since $J_l$ satisfies \ac, there exists a $J_l$-covering $R'$ of $B'$ such that the composite covering $g\circ R'$ refines $R$. $E$ is therefore locally trivial for $g\circ R'$, by \cref{lem:refinetriv}. 
Let $\bar g_{E}\colon g^*E\to E$ be a Cartesian lift of $g$ at $E$ in the fibration $P$ and $f'\colon\tilde B'\to B'$ an arrow of $R'$. Since the composite of Cartesian arrows $f'^*g^*E\to g^*E\to E$ is Cartesian over $g\circ f'$, $f'^*g^*E$ is trivial. Thus $g^*E$ is a locally trivial object.
\end{Pf}

This proposition is important because we not only get categories $\Loc_{B}$ for each locally trivial object $B$ in the base, but also functors $f^{*}\colon\Loc_{B}\to \Loc_{A}$ for each arrow $f\colon A\to B$ in $\Loc_{b}$. 

\begin{Rem}[Dualization]
  We want to give a few hints concerning the dualization of the theory. Remark first that the opposite of a fibration $P$ is an opfibration, whose opcartesian arrows are the opposite of the Cartesian arrows of $P$. Therefore, the dual of a Cartesian morphism is an opcartesian morphism. Now, an \emph{opfibred site} is a pair $(P,J)$ where $P$ is an opfibration and $J$ a \emph{co-covering function} on $P_{b}$. Together with a subfunctor of trivial objects, an opfibred site determines a subfunctor \Loc of \emph{locally trivial objects} and a co-covering function $J_{l}$ on $\Loc_{b}$. An \emph{opfibred site with trivial objects} is thus a triple $(P,J,\Triv)$ where $P$ is an opfibration, $\Triv\subset P$ is a replete and globally full sub-opfibration, $J$ is a co-covering function satisfying \am and such that $J_{l}$ satisfies \aco. A fibred site with trivial objects has the same locally trivial objects as its dual, but opposite categories: $\Loc(P^{op},J^{op},\Triv^{op})=\Loc(P,J,\Triv)^{op}$ and $(J^{op})_{l}=(J_{l})^{op}$. With this in mind, one can easily dualize the results of this section.
\end{Rem}

We now come back to \cref{ex:triv} and study locally trivial objects in these contexts. If the subfunctor is obtained by a subclass of the fibre over a terminal object and the covering function $J$ is a coverage satisfying \am, one automatically has a fibred site with trivial objects, and so we do not mention it.
\begin{Exs}[\label{ex:loctriv}]
\item\label{ex:loctriv:triv}  (\textit{Trivial cases}) Let $P$ be a fibration equipped with a subfibration of trivial objects $\Triv\subset P$. Then $(P,\mathit{Coarsest},\Triv)$ and $(P,\mathit{Finest},\Triv)$ are fibred sites with trivial objects. Moreover, the locally trivial objects are respectively the trivial objects and all objects.
\item (\textit{$G$-torsors}) Consider the context of \cref*{ex:triv}(\ref{ex:triv:gbun}). When \C is given a subcanonical pretopology, locally trivial $G$-bundles are called \emph{$G$-torsors} \cite{Vis08}. For instance, when $\C=\Top$ endowed with the pretopology of open subset coverings, $G$-torsors are precisely principal $G$-bundles. On the other hand, $G$-torsors in the pretopology of jointly surjective open maps correspond to the more general notion of principal bundles as defined in Husemoller's classical book \cite{Hus94}, \cite[Prop.\ 3.3.16]{MyThesis}.

\item (\textit{Vector bundles}) In \cref*{ex:triv}(\ref{ex:triv:vbun}), when \Top is endowed with the pretopology of open subset coverings, locally trivial objects are precisely real vector bundles. The same construction holds in the complex case.

\item (\textit{Locally constant sheaves}) When \Top has the open subset pretopology, locally trivial objects in the situation of \cref*{ex:triv}(\ref{ex:triv:csheaves}) are called \emph{locally constant sheaves}.

\item (\textit{Schemes}) Consider affine schemes in \cref*{ex:triv}(\ref{ex:triv:affine}). When \Top has the open subset pretopology, then locally trivial objects are schemes.

\item\label{ex:loctriv:proj} (\textit{Finitely generated projective modules}) Let us consider the opfibration of modules over commutative rings with its sub-opfibration over \Comm of finitely generated free modules, as in \cref*{ex:triv}(\ref{ex:triv:freemod}). We put on \Comm the Zariski Grothendieck pre-cotopology. By a classical theorem of commutative algebra \cite[Th.\ 1, n°2, §5, Ch.\ 2]{Bou98}, locally trivial objects in this situation are precisely the finitely generated projective modules \cite{MyThesis}.

\item\label{ex:loctriv:fp} (\textit{Finitely presented modules}) Consider the opfibred site $\Mod\to(\Comm,\Zar)$ and its sub-opfibration of finitely presented modules as in example \ref*{ex:triv}(\ref{ex:triv:fp}). Then, locally trivial objects coincide with the trivial ones \cite[Cor.\ of Prop.\ 3, n°1, §5, Ch.\ 2]{Bou98}.

\item\label{ex:loctriv:locfreemod} (\textit{Locally free sheaves of modules}) We deal here with the opfibration of sheaves of modules over ringed spaces and the sub-opfibration over $\Ringed^{op}$ of free modules of finite rank as defined in \cref*{ex:triv}(\ref{ex:triv:freemodringed}). We put on $\Ringed^{op}$ the Zariski pre-cotopology, that is, the Zariski pretopology on \Ringed (we denote both identically). Locally trivial objects in this context are called \emph{locally free modules of finite rank} or \emph{vector bundles} (the rank of such modules is defined locally \cite{WeiKbook}).

\item (\textit{Quasi-coherent sheaves of modules}) The opfibred site of interest is the restriction of the previous opfibred site to the (non-full) subcategory $\LRinged\subset\Ringed$ of locally ringed spaces and their morphisms. Together with the sub-opfibration of affine objects of \cref*{ex:triv}(\ref{ex:triv:affshemod}), one can show, using \cref*{lem:J_l}(\ref{lem:J_l:ii}), that this forms an opfibred site with trivial objects. Locally trivial objects in the base are affine schemes, whereas those in the total category are \emph{quasi-coherent sheaves of modules}\footnote{The category of schemes is the category of locally trivial objects in the category \LRinged of locally ringed spaces. Schemes are also the locally trivial objects in the whole category \Ringed, but the category of locally trivial objects there has too many morphisms. Restricting to the category \LRinged is in fact a shortcut that prevents us from talking about locally trivial \emph{morphisms}, which is precisely what morphisms of locally ringed spaces between schemes are. See \cite{MyThesis} for further discussion.}. If one restricts the domain of the morphism of fibrations of \cref{eq:4} to the sub-opfibration of finitely presented modules $\Mod_{c}^{fp}\to\Comm$, then one obtains \emph{coherent sheaves of modules}.\footnote{This definition agrees with Weibel \cite{WeiKbook}, but differs from \cite{Har77}, where the modules are required to be finitely generated instead of finitely presented. See \cite{WeiKbook} for a discussion on the different definitions of coherent sheaves of modules.} When considering the sub-opfibration of finitely generated free modules $\Free\to\Comm$, one recovers locally free modules of finite rank as in the preceding example, except that they now are over schemes only.

\end{Exs}

\begin{Def}\label{def:morfibsites}
  A \emph{morphism of fibred sites with trivial objects} \[(F,G)\colon (P,J,\Triv)\to (P',J',\Triv')\] is a Cartesian morphism of fibrations $(F,G)\colon P\to P'$ such that the following cube commutes
  \begin{equation}
    \label{eq:7}
    \begin{aligned}
      \xymatrix@!0@=1.2cm{
        & \Triv'_{t} \ar@{^{(}->}[rr]\ar'[d][dd] & & \E' \ar[dd]^(.3){P'}\\
        \Triv_{t} \ar[ur]\ar@{^{(}->}[rr]\ar[dd] & & \E \ar[ur]^{F}\ar[dd]^(.3){P}\\
        & \Triv'_{b} \ar@{^{(}->}'[r][rr] & & (\B',J')\\
        \Triv_{b} \ar@{^{(}->}[rr]\ar[ur] & & (\B,J) \ar[ur]_{G} }
    \end{aligned}
  \end{equation}
and $G\colon(\B,J_{l})\to(\B',J'_{l})$ is a morphism of sites.
\end{Def}

The following result is easy to prove using \cref{lem:refinetriv}.
\begin{Prop}\label{prop:morfibsites}
  A morphism of fibred sites with trivial objects as in \cref*{eq:7} induces the following commutative diagram, in which all vertical faces are Cartesian morphisms whose base functors are morphisms of sites.
  \begin{equation}
    \label{eq:8}
    \begin{aligned}
 \xymatrix@!0@=1.5cm{
&\Triv'_{t}\ar'[d][dd]\ar@{^(->}[rr] & &\Loc'_t\ar'[d][dd]\ar@{^(->}[rr]& & \E'\ar[dd]^{P'}\\
\Triv_{t}\ar[ru] \ar@{^(->}[rr] \ar[dd]& & \Loc_t\ar[ru] \ar@{^(->}[rr] \ar[dd]& & \E\ar[ru]^{F}\ar[dd]^(.4)P&\\
& (\Triv'_{b},J'_{l})\ar@{^(->}'[r][rr] & & (\Loc'_b,J'_{l})\ar@{^(->}'[r][rr] & & (\B',J'_{l})\\
(\Triv_{b},J_{l})\ar[ru] \ar@{^(->}[rr] & & (\Loc_b,J_{l})\ar[ru] \ar@{^(->}[rr] & & (\B,J_{l})\ar[ru]_{G} &
}
    \end{aligned}
\end{equation}
\cqfd
\end{Prop}

\begin{Cor}
  Let $(P,J,\Triv)$ and $(P,J',\Triv)$ be fibred sites with trivial objects such that $J_{l}\preceq J'_{l}$. Then, $\Loc(J)\subset\Loc(J')$ is a (replete, globally full) subfibration.\Par
In particular, if $J_{l}\equiv J'_{l}$, then $\Loc(J)=\Loc(J')$.
\qed
\end{Cor}

\begin{Exs}
\item One has the following morphism of opfibred sites with trivial objects.
\[\xymatrix@!0@=1.4cm{
        & \mathit{Free}^{fr} \ar@{^{(}->}[rr]\ar'[rd][rrdd]& & \Oc{}\text{-}\Mod \ar[dd]^(.3){P'}\\
        (\Free) \ar[ur]\ar@{^{(}->}[rr] \ar[rrdd]& & {\Mod_{c}} \ar[ur]\ar[dd]^(.3){P}\\
        & & & (\Ringed^{op},\mathit{Zar})\\
        & & (\Comm,\mathit{Zar}) \ar[ur]
}\]
See \cite[Lemma 2.3.7]{Qin02} for the morphism of sites in the base, and \cite[Proposition II.5.2]{Har77} for the fact that the morphism of opfibrations restricts to trivial objects. Its Cartesianness is proved in \cite{MyThesis} (part of the result is contained in the latter cited Proposition of \cite{Har77}). Therefore, there is an induced morphism on locally trivial objects. It happens to be a fibrewise equivalence of categories \cite[Ch.\ I]{WeiKbook}.
\item Pullbacks of fibred sites with trivial objects provide examples. Let $(P,J,\Triv)$ be a fibred site with trivial objects. Let $\Triv'_{b}\subset (\B',J')$ be a full subcategory of a site $(\B',J')$ such that $J'_{l}$ is a coverage. Suppose that $G\colon(\B',J'_{l})\to(\B,J_{l})$ is a morphism of sites that preserve trivial objects. Then, the following cube, whose side faces are (the canonical choices of) pullbacks, is a morphism of fibred sites with trivial objects:
\[\xymatrix@!0@=1.5cm{
        & \Triv_{t} \ar@{^{(}->}[rr]\ar'[d][dd] & & \E \ar[dd]^{P}\\
        **[gray]{\bar G^{*}\Triv_{t}} \ar@*{[gray]}[ur]\ar@*{[gray]}@{^{(}-->}[rr]\ar@*{[gray]}[dd] & & **[gray]{G^{*}\E} \ar@*{[gray]}[ur]\ar@*{[gray]}[dd]\\
        & \Triv_{b} \ar@{^{(}->}'[r][rr] & & (\B,J)\\
        \Triv'_{b} \ar@{^{(}->}[rr] \ar[ur]^{\bar G} & & (\B',J') \ar[ur]_{G}
}\]
\end{Exs}

\section{Exact structure}
\label{sec:exact-structure}
Given a fibred site with trivial objects, one would like to study locally trivial objects $B$ in the base via their associated category $\Loc_{B}$ of locally trivial objects over $B$. In fact, one would like to study the whole category $\Loc_{b}$ via the fibration $\Loc$, since a morphism $f\colon A\to B$ gives rise to a functor $f^{*}\colon\Loc_{B}\to\Loc_{A}$. Moreover, if the categories $\Loc_{B}$ (and some of the inverse image functors between them) have the required properties and structures, one can apply to them a $K$-theory machine, transferring the study into the realm of infinite loop spaces. So now, we must find conditions so that the tools introduced so far provide the properties and structure required by the $K$-theory functor. 

Firstly, the category $\Loc_{B}$ of locally trivial objects over $B\in \Loc_{b}$ must be skeletally small, and this will be an hypothesis. Secondly, $\Loc_B$ should come with some suitable structure.  Here, we study conditions that ensure an exact structure on the categories of locally trivial objects.\footnote{In some examples, such as ringed spectra, one has to work with a general Waldhausen structure that does not come from an exact one. Moreover, the definition of the $K$-theory of a space via its category of (real or complex) vector bundles requires taking not only the exact structure, but also the topological enrichment, into account (see \cite[Chapter 5]{MyThesis} and references therein for more details). We do not study these situations in this article.} Finally, we turn to conditions ensuring exactness of inverse image functors.

\subsection{Preliminaries about exact categories}
\label{sec:prel-about-exact}
We start with some preliminary material about exact categories, in order to fix the notation and terminology. This is a very brief account. Among good references are \cite{BKII00,BorII94,Osb00,Qui73}.
\begin{Defs}
  \item A \emph{preadditive category} is a category enriched over the monoidal category \Ab of abelian groups. An \emph{additive functor} is an \Ab-enriched functor.
  \item A \emph{preadditive subcategory} of a preadditive category \C is a subcategory $\A\subset \C$ with a preadditive structure and whose inclusion functor is additive.
  \item An \emph{additive category} is a preadditive category with a zero object and with a biproduct of any pair of objects.
  \item\label{def:addsubcat} An \emph{additive subcategory} of an additive category \C is a preadditive subcategory ${\A\subset\C}$ that is additive.
\end{Defs}

\begin{Rem}
An additive functor automatically preserves zero objects and biproducts. Therefore, the zero objects and biproducts of an additive subcategory $\A\subset\C$ are necessarily zero objects and biproducts in the ambient category \C. If moreover the subcategory \A is replete, then it is closed under taking zero objects and biproducts of objects of \A in \C.
\end{Rem}

\begin{Defs}
\item A \emph{short exact sequence} in an preadditive category \C is a sequence
  \begin{equation}
    \label{eq:9}
    A\ra f B\ra g C
  \end{equation}
of morphisms in \C such that $f$ is the kernel of $g$ and $g$ the cokernel of $f$. It is \emph{split} if either $f$ is a split mono or $g$ is a split epi.
\item A \emph{pre-exact category}\footnote{This is called a $G$-exact category in \cite{BKII00}.} is a pair $(\C,\se)$ where \C is an additive category and Ex is a chosen class of short exact sequences in \C, called \emph{admissible short exact sequences}, such that:
  \begin{enumerate}
  \item All split exact sequences belong to \se,
  \item \se is closed under isomorphisms: given a short exact sequence \eqref{eq:9} of \se and a commutative diagram whose vertical arrows are isomorphisms
\[
\xymatrix{A\ar[r]^{f}\ar[d]_{\cong} & B\ar[r]^{g}\ar[d]_{\cong} & C\ar[d]^{\cong}\\
A'\ar[r] & B'\ar[r] & C'
}
\]
the bottom sequence (which is automatically short exact) belongs to \se.
\item \se is closed under direct sums: if $A\ra f B\ra g C$ and $A'\ra{f'} B'\ra{g'} C'$ are in \se, then their direct sum (which is automatically short exact) is in \se.
  \end{enumerate}
\item An \emph{exact functor} \[F\colon(\A,\se_{\A})\to(\B,\se_{\B})\] between pre-exact categories is an additive functor $F\colon\A\to\B$ such that $$F(\se_{\A})\subset \se_{\B},$$ that is, it sends a short exact sequence of $\se_{\A}$ to a short exact sequence in \B that moreover belongs to $\se_{\B}$.
\item A \emph{pre-exact subcategory} of a pre-exact category $(\C,\se_{\C})$ is a pre-exact category $(\A,\se_{\A})$ such that $\A\subset\C$ is an additive subcategory and the inclusion functor is exact.
\end{Defs}

\begin{Rem}
If a short exact sequence $A\ra f B\ra g C$ splits, then both morphisms $f$ and $g$ admit a splitting, and $B$ is a biproduct of $A$ and $C$.
\end{Rem}

\begin{Exs}\label{ex:preex}
  \item Any additive category with its class of all split short exact sequences is pre-exact and called \emph{split pre-exact}. On the other hand, any additive category together with its class of all short exact sequences is a pre-exact category, called \emph{repletely pre-exact category}. When no class of short exact sequence is specified, then we mean the repletely pre-exact structure.
  \item\label{ex:preex:2} A full additive subcategory \A of a pre-exact category $(\C,\se)$ together with the class of all short exact sequences of \se whose objects belong to \A is a pre-exact category.
\end{Exs}

We give, for the sake of simplicity, a definition of an exact category via an embedding in an abelian category. There is a purely axiomatic way of defining an exact category (see for instance \cite{BKII00}), but it is quite long to state. It is slightly more general in the sense that all exact categories as defined below satisfy these axioms, but only skeletally small ``axiomatically defined'' exact categories are exact categories in the sense given here (via a sort of Yoneda embedding).
\begin{Def}\label{def:exact}
 An \emph{exact category} (in the sense of Quillen) is a pre-exact category $(\C,\se)$ such that there exists an abelian category \D for which the following holds:
   \begin{enumerate}[(i)]
   \item \C is a full additive subcategory of \D.
   \item \se is the class of all short exact sequences in \D whose objects belong to \C.
   \item \label{def:exact:closed} \C is closed under extensions in \D: Given a short exact sequence
\[A\ra f B\ra g C\]
in \D such that $A$ and $C$ are isomorphic to objects of \C, then $B$ must be isomorphic to an object of \C.
   \end{enumerate}
\end{Def}

\begin{Rems}[\label{rem:exact}]
\item Any skeletally small split pre-exact category is exact (via some variant of the Yoneda embedding, see \cite{Qui73}), and thus is called \emph{split exact}.
  \item Any abelian category with its repletely pre-exact structure is trivially exact. \emph{By convention, an abelian category is considered as an exact category when equipped  with this structure.}
\item\label{rem:exact:subcat} Let $(\C,\se_{\C})$ be an exact category. Let $\B\subset\C$ be a full additive subcategory and let us equip it with the class $\se_{\B}$ of all short exact sequences of $\se_{\C}$ whose objects are in \B. Suppose that \B is closed under extensions in $\se_{\C}$. This is a pre-exact category, as we have already noticed, and thus an exact category. Following the literature, we call this an \emph{exact subcategory} of the exact category $(\C,\se_{\C})$. An exact category is always an exact subcategory of the abelian category that determines its structure.
\end{Rems}

The definition of an exact category is somehow redundant. Here is minimal characterization that is useful in practice.

\begin{Lem}\label{lem:exactcat}
  Let \D be an abelian category and $\C\subset\D$ be a full subcategory containing a zero object of \D. Suppose that \C is closed under extensions in \D. Then, \C is an additive subcategory of \D and an exact category with respect to short exact sequences in \D with objects in \C.
\end{Lem}
\begin{Pf}
Note that a full subcategory \C of a preadditive category \D is automatically a preadditive subcategory. Moreover, the closure of \C under extensions in \D implies
  that \C has biproducts of all pairs of objects. Finally, by \cref{ex:preex}, any full additive subcategory \C of an abelian
  category \D is a pre-exact category with respect to exact sequences
  of \D belonging to \C.
\end{Pf}
  Thus, there no need to check the pre-exact
  axioms. Still, one should remember that, whereas it is a property of a pre-exact category, exactness is a structure of the underlying category.

\subsection{Exactness of locally trivial objects}
\label{sec:exactn-locally-triv}
Consider a fibred site with trivial objects $(P,J,\Triv)$. Suppose that $P$ is fibrewise abelian with inverse image functors additive. We would like its subfibration $\Loc$ of locally trivial objects to be fibrewise exact, so that one can apply $K$-theory to its fibres. The strategy for obtaining this property is the following: if the property holds for trivial objects, then it is, under some conditions, also true for locally trivial objects.

\begin{Defs}\label{def:abfibsite}
\item An \emph{additive fibration} is a category fibred in additive categories, that is, a fibration whose fibres and inverse image functors are additive. It is \emph{abelian} if moreover each fibre is abelian. 

\item An \emph{additive subfibration} of an additive fibration is a \hyperref[def:subfib]{subfibration} whose fibres are \hyperref[def:addsubcat]{additive subcategories} (it is then automatically additive itself). 

\item An \emph{exact fibration} $P$ is an additive fibration whose fibres are endowed with an exact structure. We do not require in general that its inverse image functors be exact.
\item A morphism $f\colon A\to B$ in $P_{b}$ is \emph{$P$-flat} or \emph{flat in $P$} if its inverse image functor is exact. A covering $R$ of an object of $P_{b}$ is \emph{$P$-flat} or \emph{flat in $P$} if all its arrows are $P$-flat.

\item An \emph{exact subfibration} $Q$ of an exact fibration $P$ is a \hyperref[def:subfib]{subfibration} whose fibres are \hyperref[rem:exact:subcat]{exact subcategories} (see remark \ref*{rem:exact}(\ref{rem:exact:subcat})).

\item\label{def:abfibsite:triv} An \emph{abelian fibred site with trivial objects} is a \hyperref[def:fibsitetriv]{fibred site with trivial objects} \linebreak $(P,J,\Triv)$ such that: 
  \begin{enumerate}
  \item $P$ is an abelian fibration.
  \item $\Triv\subset P$ is an exact subfibration of the abelian fibration $P$.
  \item $J_{l}$-coverings are flat in $P$. 
  \item\label{def:abfibsite:triv:3} The covering function $J_{l}$ satisfies (L).
  \end{enumerate}
\end{Defs}

We now state the main result of this section.

\begin{Prop}\label{prop:abfibsite}
  Consider an abelian fibred site  with trivial objects $(P\colon\E\to\B,J,\Triv)$.\Par
  Then, the subfibration of locally trivial objects is an exact subfibration of $P$.
\end{Prop}

\begin{Pf}
We already know that locally trivial objects form a replete and globally full subfibration of $P$ (\cref{prop:loc}). In particular, for each $B\in\Loc_{b}$, one has a full subcategory $\Loc_{B}\subset\E_{B}$. By \cref{lem:exactcat}, one only needs to prove that $\Loc_{B}$ contains a zero object of $\E_{B}$ and is closed under extensions. $\Loc_{B}$ contains any zero object $0$ of $\E_{B}$. Indeed, let $R$ be a $J_{l}$-covering of $B$ and $f\colon A\to B$ in $R$. Inverse image functors being additive, they preserve zero objects. Thus, $f^{*}0$ is a zero object of $\E_{A}$. Again, $\Triv_{A}\subset\E_{A}$ being a replete additive subcategory, it contains $f^{*}0$. Thus, $0$ is locally trivial. 

It remains to show the closure under extensions. Consider a short exact sequence $E'\to E\to E''$ in the abelian category $\E_{B}$ where $E'$ and $E''$ are locally trivial (recall that $\Loc_{B}$ is replete). Let $R$ be a $J$-covering that trivializes both $E'$ and $E''$, which exists by \cref{lem:comtriv}. Let $f\colon A\to B$ be in $R$. By hypothesis, the inverse image functor $f^{*}$ is exact and thus, one obtains a short exact sequence
\[f^{*}E'\to f^{*}E\to f^{*}E''\]
in $\E_{A}$. Now, since $f^{*}E'$ and $f^{*}E''$ are trivial and since $\Triv_{A}\subset\E_{A}$ is an exact subcategory, $f^{*}E$ is also trivial. Therefore, $E$ is locally trivial.
\end{Pf}

\begin{Rem}[Dualization]
 The straightforward dualization of the theory in this section relies on the fact that, given an exact category $\C$ defined by its inclusion in an abelian category \A, its dual $\C^{op}$ is still an exact category, defined by its inclusion the abelian category $\A^{op}$. 
\end{Rem}

Recall that in many cases, the opfibred site comes from the opfibration of modules over monoids in a monoidal opfibration. Before studying applications of \cref{prop:abfibsite}, we determine what is needed in order to obtain an abelian structure on this opfibration.

\begin{Thm}[Ardizzoni, \cite{Ard04}]\label{thm:ardi}
    Let \V be a monoidal category whose underlying category is abelian. Let $R$ be a monoid in \V. Suppose that the functor $-\otimes R$ preserves finite colimits. \Par
Then, the category $\Mod_{R}$ of $R$-modules in \V is abelian.
  \end{Thm}
  \begin{Def}
    A (right) \emph{monoidal abelian opfibration} is a monoidal opfibration whose underlying opfibration is abelian and such that each tensor
    $-\otimes E$ is additive.
  \end{Def}
  \begin{Cor}\label{cor:monabopfib}
    Let $\E\to\B$ be a monoidal abelian opfibration. Suppose that the direct image functors and the right tensors preserve cokernels.\Par
Then, the opfibration of modules $\Mod(P)\to\Mon(P)$ is abelian.
\Par
The same is true in the commutative setting.
  \end{Cor}
  \begin{Pf}
    The fact that its fibres are abelian is a direct consequence of \cref{prop:modovermon} and \cref{thm:ardi}. Let $\phi\colon R\to S$ be a morphism of monoids in $P$. Let us write $f\coloneqq P(\phi)$. The direct image functor associated to $\phi$ is, following the proof of \cref{prop:modovermon}, given by the composite:
\[\Mod_{R}\ra{f_{*}}\Mod_{f_{*}R}\ra{-\otimes_{f_{*}R}S}\Mod_{S}.\] The first functor is additive because $P$ is abelian. The second functor is additive, because it is a left adjoint.
  \end{Pf}

We now illustrate these result by examples. Sometimes, there is a direct proof of the exactness of the subfibration of locally trivial objects that is much simpler. Yet, these examples still show the legitimacy of the \cref{prop:abfibsite}'s hypotheses. In particular, it should become clear that requiring inverse image functors over arrows of $J_{l}$-coverings to be exact is meaningful, whereas requiring all inverse image functors to be so would be wrong.
\begin{Exs}\label{exfibsiteexact}
\item \label{exfibsiteexact:triv} (\textit{Trivial cases}) Let $\Triv\subset P$ be a subfibration of trivial objects. Suppose that $P$ is abelian and $\Triv$ an exact subfibration of $P$. Then $(P,\mathit{Coarsest},\Triv)$ and $(P,\mathit{Finest},\Triv)$ are abelian fibred sites with trivial objects (see example \ref{ex:loctriv}(\ref{ex:loctriv:triv})).
\item\label{ex:modcommexact} (\textit{Modules over commutative rings})
This is a dual example. Consider the opfibred site of modules over commutative rings with the Zariski pre-cotopology. Let us consider the sub-opfibration of trivial objects given by all finitely generated free modules over commutative rings (example \ref*{ex:loctriv}(\ref{ex:loctriv:proj})).
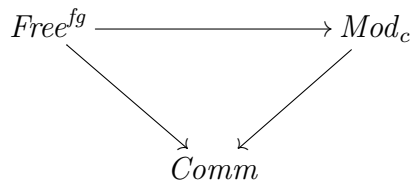
\[\xymatrix@R=1.3cm{\Free\ar[rr]\ar[rd] && \Mod_{c}\ar[ld]\\
&\Comm&
}\] 
We check that this is an \hyperref[def:abfibsite:triv]{abelian opfibred site with trivial objects}.

The opfibration $\Mod_c$ is abelian, since $\Ab\to\one$ satisfies the hypotheses of \cref{cor:monabopfib}. Moreover, each category $\Free_{R}$ of finitely generated free $R$-modules is an exact subcategory of the abelian category $\Mod_{R}$ of all $R$-modules, since free modules are projective. 

Now, given an element $a\in R$, the localization morphism $R\to R[a^{-1}]$ is flat, that is, the corresponding extension of scalars functor is exact \cite{Bou98}. Note that it is not true that all extension of scalars functors are flat. Finally, since all commutative rings are considered as trivial in the base, $J=J_{l}$ in this case. Since $J$ is a pre-cotopology, it satisfies \alo. In conclusion, $(\Mod_{c},\Zar,\Free)$ is an abelian opfibred site with trivial objects.

When extending the trivial objects to all finitely presented modules, one still obtains an abelian opfibred site with trivial objects $(\Mod_{c},\Zar,\Mod_{c}^{\mathit{fp}})$. Indeed, finitely presented modules over a commutative ring $R$ form an exact subcategory in the category of modules over $R$ \cite[Lemma 7.5.4]{StaProj}.

\item\label{ex:modringedexact} (\textit{Modules over ringed spaces})
Consider now the opfibred site of sheaves of modules over ringed spaces in the Zariski pre-cotopology $\Oc{}$-$\Mod\to(\Ringed^{op},\Zar)$, equipped with the sub-opfibration over $\Ringed^{op}$ of free modules of finite rank (example \ref{ex:loctriv}(\ref{ex:loctriv:locfreemod})). In this case, $J_{l}=J$ and thus all conditions relative to the site are satisfied. This opfibration is abelian, by \cref{cor:monabopfib}. Indeed, the opfibration $\Sh\to\Top^{op}$ is a bifibration whose fibres are closed monoidal abelian categories. It is an exact sub-opfibration of the abelian opfibration of all modules (because exactness can be checked stalkwise).

Now, recall that a morphism of ringed spaces $(f,f^{\sharp})\colon(X,\Oc X)\to (Y,\Oc Y)$ is said to be \emph{flat} if for each $x\in X$, the canonical homomorphism of rings
$\Oc{Y,f(x)}\to\Oc{X,x}$ is flat. Such a morphism is also flat in the sense of \cref{def:abfibsite} (see \cite[Lemma 13.17.2]{StaProj}). Since open immersions are flat morphisms of ringed spaces \cite{Qin02}, Zariski coverings are flat. Therefore, $(\Oc{}$-$\Mod,\Zar,\mathit{Free}^{\mathit{fr}})$ is an abelian opfibred site with trivial objects.

\item\label{ex:modschexact} (\textit{Modules over schemes}) 
We now treat the situation of the opfibred site with trivial objects $(\Oc{}\text{-}\Mod_{l},\Zar,\Triv)$, where $\Triv$ is induced by the opfibrations $\Mod_{c}$, $\Mod_{c}^{\mathit{fp}}$ or $\Free$. We check that these form abelian opfibred sites with trivial objects. Since, in order to prove it, we are using the fact that given a scheme $(X,\Oc X)$, quasicoherent sheaves of $\Oc X$-modules form an exact subcategory of the category of all sheaves of $\Oc X$-modules, we do not claim to give a new proof of this exactness result using \cref{prop:abfibsite}! However, we use \cref{prop:abfibsite} to provide a proof that coherent sheaves of modules (as defined here) form an exact subcategory, and we are not aware of another proof of this fact. 

The opfibration $\Oc{}\text{-}\Mod_{l}$ is abelian as a sub-opfibration of $\Oc{}\text-\Mod$. Moreover, coverings of $\Zar_{l}$ are flat, since all Zariski coverings are so. In addition, $\Zar_{l}$ satisfies $(L)$ by \cref*{lem:J_l}(\ref{lem:J_l:ii}). It remains to show that \Triv is an exact sub-opfibration of $\Oc{}$-$\Mod_{l}$.

Let $(X,\Oc X)$ be an affine scheme. Let $0\to\mathcal{F}'\to\mathcal{F}\to\mathcal{F}''\to 0$ be a short exact sequence of $\Oc X$-modules with $\mathcal{F}'$ and $\mathcal{F}''$ affine. Quasicoherent sheaves of $\Oc X$-modules form an exact subcategory of the category of all sheaves of $\Oc X$-modules \cite[Ch.\ 5, Prop.\ 1.12(d)]{Qin02}\footnote{Liu's definition of quasicoherent modules in \cite{Qin02} is equivalent to ours over schemes (see \cite[Prop.\ 5.4]{Har77} and \cite[Ch.\ 5, Th.\ 1.7]{Qin02}).} and therefore, $\mathcal{F}$ is quasicoherent. Moreover, an $\Oc Y$-module over a scheme $(Y,\Oc Y)$ is quasicoherent if and only it is locally affine over all Zariski coverings of $(Y,\Oc Y)$ by affine subschemes \cite[Ch.\ 5, Th.\ 1.7]{Qin02}. Therefore, $\mathcal{F}$ is affine.

Suppose now that $\mathcal{F}'$ and $\mathcal{F}''$ are finitely presented affine modules. Since over affine schemes, the global section functor preserves short exact sequences with left term quasicoherent \cite[Ch.\ 5, Prop.\ 1.8]{Qin02}, the sequence \[0\to\mathcal{F}'(X)\to\mathcal{F}(X)\to\mathcal{F}''(X)\to 0\] is an exact sequence of  $\Oc X(X)$-modules. Now, finitely presented modules form an exact subcategory in the category of modules (see \cref*{exfibsiteexact}(\ref{ex:modcommexact})). Therefore, $\mathcal{F}(X)$ is finitely presented. Since $\mathcal{F}$ is quasicoherent, \cite[Ch.\ 5, Th.\ 1.7]{Qin02} applies to show that $\mathcal{F}$ is a finitely presented affine object.

If \Triv is induced by the opfibration $\Free$, then \Triv is the restriction of $\mathit{Free}^{\mathit{fr}}$ to the category of affine schemes. It is an exact sub-opfibration by the preceding example.

Note that one can extend the Grothendieck pretopologies of \cref*{ex:pretop}(\ref{ex:pretop:sch}) to pretopologies satisfying \tc and \tl on \LRinged by assigning just the identity covering to a non-schematic locally ringed space. Each of these pretopologies gives rise to an abelian opfibred site with trivial objects $(\Oc{}\text{-}\Mod,J,\Triv)$ when \Triv is induced by $\Mod_{c}$, $\Mod_{c}^{\mathit{fp}}$ or $\Free$. Indeed, \cref{lem:J_l}(\ref{lem:J_l:ii}) applies to all of them, and all of their morphisms are flat \cite[Ch.\ 24 and 30]{StaProj}.

\end{Exs}

\subsection{Exactness of inverse image functors}
\label{sec:exactn-inverse-image}
Under the conditions of \cref{prop:abfibsite}, categories of locally trivial objects are exact categories. \Cref{prop:abfibsite} also guarantees that inverse image functors are additive. Yet, in order to induce a morphism in $K$-theory, a functor must be exact. Thus, only flat morphisms of locally trivial objects induce morphisms in $K$-theory. We show here that under some hypotheses, the flatness of morphisms of trivial objects implies flatness of morphisms of locally trivial objects.

\begin{Defs}
[Let $P$ be an exact fibration.]
\item A morphism $f$ in $P_{b}$ is \emph{faithfully flat in $P$} if its inverse image functor $f^{*}$ preserves and reflects short exact sequences.
\item A covering $R$ of an object of $P_{b}$ is \emph{faithfully flat in $(P,J)$} if its cartesian lift ${R^{*}=\{f^{*}\mid f\in R\}}$ collectively preserves and reflects exact sequences.
\end{Defs}

\begin{Prop}\label{prop:exactinvimage}
  Let $(P,J,\Triv)$ be an abelian fibred site with trivial objects and $\Loc\subset P$ be its exact subfibration of locally trivial objects. Suppose that $J_{l}$-coverings are faithfully flat in \Loc and that morphisms of trivial objects are flat in \Triv.
\Par
Then, morphisms of locally trivial objects are flat in \Loc.
\end{Prop}
\begin{Pf}
Let $f\colon A\to B$ be a morphism of locally trivial objects. Let 
\[s=\{0\to E'\to E\to E''\to 0\}\]
be an exact sequence in $\Loc_{B}$. By \cref{lem:comtriv}, there is a $J$-covering 
\[S=\{h_{j}\colon B_{j}\to B\}_{j\in J}\] 
that trivializes $E''$, $E'$ and $E$ together. Since $J_{l}$ satisfies \ac, there is a $J_{l}$-covering ${R=\{g_{i}\colon A_{i}\to A\}_{i\in I}}$ of $A$ such that, for all index $i\in I$, there is an index $j\in J$ and a morphism $f_{i}$ in \Triv that make the following square commute:
\[\xymatrix@=1.5cm{
A_{i}\ar@{-->}[r]^{f_{i}}\ar[d]_{g_{i}} & B_{j}\ar[d]^{h_{j}}\\
A\ar[r]_{f} & B.
}\]
Thus, for all $i\in I$, $g_{i}^{*}\circ f^{*}\cong f_{i}^{*}\circ h_{j}^{*}$. The covering $R$ is $P$-flat since it belongs to $J_{l}$. Moreover, each $f_{i}$ is flat in \Triv by hypothesis. Therefore, $g_{i}^{*}\circ f^{*}(s)$ is an exact sequence in $\Loc_{A_{i}}$ for all $i\in I$. Since $J_{l}$-coverings are supposed to be faithfully flat in \Loc, $f^{*}(s)$ is exact.
\end{Pf}
\begin{Exs}
\label{ex:loctrivflat}
\item Consider the abelian opfibred site with trivial objects \[(\Mod_{c},\Zar,\Free).\] It is a trivial fact that all homomorphisms of rings are flat in \Loc, since the subcategories of finitely generated projective modules are split exact. Nonetheless, note that \cref{prop:exactinvimage} applies. Indeed, Zariski coverings are faithfully flat by \cite[Prop.\ 3, n°1, §5, Ch.\ 2]{Bou98}.
\item\label{ex:loctrivflatsheaves} We turn now to the abelian opfibred site with trivial objects $(\Oc{}$-$\Mod,\Zar,\mathit{Free^{\mathit{fr}}}$). The categories of locally trivial objects in this context, which are locally free sheaves, are not always split exact \cite[Ex.\ 7.1.3, Ch.\ II]{WeiKbook}. \Cref{prop:exactinvimage} applies to show that all morphisms of ringed spaces are flat in \Loc. 
\end{Exs}

\section{\texorpdfstring{$K$-theory}{K-theory}}
\label{sec:k-theory}
We are now ready to define the $K$-theory of a locally trivial object in an abelian fibred site with trivial objects.
\subsection{Definition}
\label{ssec:definitions}

In this part, $(P,J,\Triv)$ is an abelian (op-)fibred site with trivial objects and $\Loc\subset P$ denotes the exact sub-(op-)fibration of locally trivial objects. \emph{We suppose that $\Loc$ is fibrewise skeletally small.}

\begin{Def}
The \emph{$K$-theory of a locally trivial object $B\in\Loc_{b}$} is the $K$-theory of the exact category $\Loc_{B}$: \[K_{(P,J,\Triv)}(B)=K(B)\coloneqq K(\Loc_{B}).\]
\end{Def}

\begin{Rem}
In important examples, the abelian opfibration $P$ is built from a module opfibration. Consider a monoidal abelian opfibration $Q$ satisfying the conditions of \cref{cor:monabopfib}. Then, the opfibration of modules $\Mod(Q)$ is abelian. One obtains various abelian opfibrations $P$ by pullbacks of $\Mod(Q)$ along functors $F\colon\B\to Q_{b}$.
\end{Rem}

\subsection{Morphisms in \texorpdfstring{$K$-theory}{K-theory}}
\label{ssec:morphisms-k-theory}
We now investigate the question of morphisms in $K$-theory induced by such contexts. Consider a fibred situation (the dual works similarly). Since $\Loc$ is an additive sub-fibration of $P$, each morphism of locally trivial objects $f\colon A\to B$ induces an additive inverse image functor $f^{*}\colon\Loc_{B}\to\Loc_{A}$. When $f$ is flat\footnote{Recall that under hypotheses of \cref{prop:exactinvimage}, all morphisms of locally trivial objects are flat.}, one thus obtains a morphism in $K$-theory
\[f^{*}\colon K(B)\to K(A).\]

We want now to compare the $K$-theories induced by two different abelian fibred sites with trivial objects.

\begin{Defs}
\item A morphism of additive (\resp exact) fibrations $(F,G)\colon P\to P'$ is \emph{additive} (\resp \emph{exact}) if it is fibrewise additive (\resp fibrewise exact).
\item  A \emph{morphism of abelian fibred sites with trivial objects} is a \hyperref[def:morfibsites]{morphism of fibred sites with trivial objects} that is additive.
\end{Defs}

Let $(F,G)\colon (P,J,\Triv)\to (P',J',\Triv')$ be a morphism of abelian fibred sites with trivial objects. Suppose that the induced morphism $\Loc\to\Loc'$ given by \cref{prop:morfibsites} is exact (for instance, if $(F,G)$ is itself exact). Then, one has, for each locally trivial object $B\in\Loc_{b}$, a morphism
\[F\colon K(B)\to K(G(B)).\]

We consider two particular cases of this situation. In the first case, the sole covering function varies. Let $(P,J,\Triv)$ and $(P,J',\Triv)$ be abelian fibred sites with trivial objects such that $J_{l}\preceq J'_{l}$ (since only the covering function varies, we omit the rest of the data in the sequel). Then, $\Loc(J)\subset\Loc(J')$ is an exact subfibration. Therefore, one has the following morphisms in $K$-theory, for each $B\in\Loc(J)_{b}$: 
\[K_{J}(B)\to K_{J'}(B).\]
Now, recall that for every abelian fibred site with trivial objects $(P, J, \Triv)$, the coarsest and finest pretopologies determine, as locally trivial objects, respectively the trivial objects and all objects (\cref{ex:loctriv}(\ref{ex:loctriv:triv}) and \ref{exfibsiteexact}(\ref{exfibsiteexact:triv})). Since $\mathit{Coarsest}_{l}\preceq  J_{l}$, there is a morphism from the trivial $K$-theory $K^{\mathit{triv}}(B)\coloneqq K_{\mathit{Coarsest}}(B)$:
\[K^{\mathit{triv}}(B)\to K(B).\]
When $P$ itself is skeletally small, since $J_{l}\preceq\mathit{Finest}_{l}$, there is a morphism into the ``total'' $K$-theory $K^{total}(B)\coloneqq K_{\mathit{Finest}}(B)$:
\[K(B)\to K^{\mathit{total}}(B).\]
In the second case, the sole subfibration of trivial objects varies.  Let $(P,J,\Triv)$ and $(P,J,\Triv')$ be abelian fibred sites with trivial objects such that $\Triv\subset \Triv'$. Then, $\Loc(\Triv)\subset \Loc(\Triv')$ is an exact subfibration. Therefore, one has the following morphisms in $K$-theory (we omit in the subscript the non-varying part):
\[K_{\Triv}(B)\to K_{\Triv'}(B).\]

\subsection{Examples}
\label{ssec:examples}
\subsubsection{ \texorpdfstring{$K$-theories of commutative rings}{K-theories of commutative rings}}
\label{sssec:k-theor-comm}

We consider here the context $(\Mod_{c},\Zar,\Free)$ of \cref*{exfibsiteexact}(\ref{ex:modcommexact}). It is the abelian opfibred site of modules over commutative rings with the Zariski pre-cotopology and with trivial objects given by finitely generated free modules. This context provides the usual $K$-theory of a commutative ring:
\[K(R)=K_{(\Mod_{c},\Zar,\Free)}(R).\]
The free $K$-theory arises from $(\Mod_{c},\mathit{Coarsest},\Free)$ and thus the discussion of \cref{ssec:morphisms-k-theory} provides the well-known morphism 
\[K^{\mathit{free}}(R)\to K(R).\] 
When considering the context $(\Mod_{c},\Zar,\Mod_{c}^{\mathit{fp}})$ where trivial objects are all finitely presented modules, locally trivial objects are the trivial ones. This gives rise to the finitely presented $K$-theory $K^{\mathit{fp}}(R)$ of the ring $R$. The discussion in \cref{ssec:morphisms-k-theory} provides us with a \emph{Cartan} morphism \[K(R)\to K^{\mathit{fp}}(R).\]

\subsubsection{ \texorpdfstring{$K$-theories of ringed spaces}{K-theories of ringed spaces}}
\label{sssec:k-theories-ringed}
We consider the abelian opfibred site with trivial objects $(\Oc{}$-$\Mod,\Zar,\mathit{Free}^{\mathit{fr}})$ of \cref*{exfibsiteexact}(\ref{ex:modringedexact}). This determines the usual $K$-theory of a ringed space:
\[K(X,\Oc X)=K_{(\Oc{}\text{-}\Mod,\Zar,\mathit{Free}^{\mathit{fr}})}(X,\Oc X).\]
Any morphism of ringed spaces $f\colon(X,\Oc X)\to (Y,\Oc Y)$ induces a morphism in $K$-theory, by example \ref*{ex:loctrivflat}(\ref{ex:loctrivflatsheaves}):
\[f^{*}\colon K(Y,\Oc Y)\to K(X,\Oc X).\]
Following \cref{ssec:morphisms-k-theory}, the coarsest pretopology determines a free $K$-theory $K^{\mathit{free}}(X,\Oc X)$ and a morphism
\[K^{\mathit{free}}(X,\Oc X)\to K(X,\Oc X).\] 

\subsubsection{ \texorpdfstring{$K$-theories of schemes}{K-theories of schemes}}
\label{sssec:k-theories-schemes}

The standard $K$-theory of a scheme $(X,\Oc X)$ is defined via its category of locally free sheaves of modules, as for general ringed spaces. It therefore comes from the abelian opfibred site with trivial objects $(\Oc{}\text{-}\Mod_{l},\Zar,\Free)$, where, by abuse of notation, we have identified $\Free$ with the sub-opfibration of $\Oc{}$-$\Mod_{l}$ that it induces (see \cref*{exfibsiteexact}(\ref{ex:modschexact})).
\[K(X,\Oc X)=K_{(\Oc{}\text{-}\Mod_{l},\Zar,\Free)}(X,\Oc X).\]

Similarly to the finitely presented $K$-theory of a commutative ring, there is a \emph{coherent} $K$-theory $K^{c}(X,\Oc X)$ of a scheme $(X,\Oc X)$, that is defined via its category of coherent modules. It therefore derives from the abelian fibred site with trivial objects $(\Oc{}\text{-}\Mod_{l},\Zar,\Mod_{c}^{\mathit{fp}})$. Following \cref{ssec:morphisms-k-theory}, one thus obtains a \emph{Cartan} morphism:
\[K(X,\Oc X)\to K^{c}(X,\Oc X).\]

Instead of modifying the trivial objects, one can change the site. But we now have a wide range of available pre-topologies $J$ on \LRinged, all giving rise to an abelian site with trivial objects with both the kinds of trivial objects mentioned above. Thus, for each pretopology $J$ on \LRinged amongst
\[\mathit{Coarsest}\subset\mathit{Zar}\subset\mathit{\acute etale}\subset\mathit{smooth}\subset\mathit{syntomic}\subset\mathit{fppf}\subset\mathit{fpqc},\]
there are $K$-theories 
\[K_{J}(X,\Oc X)=K_{(\Oc{}\text{-}\Mod_{l},J,\Free})(X,\Oc X)\]
 and 
\[K_{J}^{c}(X,\Oc X)=K_{(\Oc{}\text{-}\Mod_{l},\Zar,\Mod_{c}^{\mathit{fp}}})(X,\Oc X).\]
When $J=\mathit{Coarsest}$, we denote $K^{\mathit{free}}\coloneqq K_{\mathit{Coarsest}}$ and $K^{\mathit{aff}}\coloneqq K^{c}_{\mathit{Coarsest}}$. For each scheme $(X,\Oc X)$, one has the following commutative diagram in $K$-theory, where $X$ abbreviates $(X,\Oc X)$:
\[\xymatrix@C=.55cm@R=1.5cm{
K^{\mathit{free}}(X)\ar[r]\ar[d] & K(X)\ar[r]\ar[d] & K_{\mathit{\acute{e}t}}(X)\ar[r]\ar[d]& K_{\mathit{smth}}(X)\ar[r]\ar[d] & K_{\mathit{synt}}(X)\ar[r]\ar[d] & K_{fppf}(X)\ar[r]\ar[d] & K_{fpqc}(X)\ar[d]\\
K^{\mathit{aff}}(X)\ar[r] & K^{c}(X)\ar[r] & K^{c}_{\mathit{\acute{e}t}}(X)\ar[r]& K^{c}_{\mathit{smth}}(X)\ar[r] & K^{c}_{\mathit{synt}}(X)\ar[r] & K^{c}_{fppf}(X)\ar[r] & K^{c}_{fpqc}(X).
}\]

\addcontentsline{toc}{section}{Bibliography}
\bibliography{/Users/nimichel/Documents/Dropbox/Recherche_maths/Bibliography/biblio}

\end{document}